\documentclass[11pt,reqno]{amsart}
\usepackage[T1]{fontenc}
\usepackage[a4paper,margin=28mm]{geometry}
\usepackage{amsmath,amssymb,mathtools,amsthm}
\usepackage{microtype}
\usepackage{xcolor}
\usepackage[hidelinks]{hyperref}

\usepackage{esint}
\usepackage{enumitem}

\newcommand{\R}{\mathbb R}
\newcommand{\Rp}{\R^n_+}
\newcommand{\dd}{\,\mathrm d}
\newcommand{\loc}{\mathrm{loc}}
\newcommand{\D}{\mathcal D}
\newcommand{\Sch}{\mathcal S}
\newcommand{\Pclass}{\mathcal P}
\newcommand{\bu}{\mathbf u}
\newcommand{\bU}{\mathbf U}

\newcommand{\bw}{\mathbf w}
\newcommand{\bF}{\mathbf F}
\newcommand{\bC}{\mathsf C}
\newcommand{\ind}{\mathbf 1}
\newtheorem{theorem}{Theorem}[section]
\newtheorem{proposition}[theorem]{Proposition}
\newtheorem{lemma}[theorem]{Lemma}

\theoremstyle{definition}
\newtheorem{definition}[theorem]{Definition}
\theoremstyle{remark}
\newtheorem{remark}[theorem]{Remark}
\numberwithin{equation}{section}

\hypersetup{colorlinks=true,linkcolor=blue,urlcolor=blue}
\newcommand{\be}{\begin{equation} \label}
\newcommand{\ee}{\end{equation}}

\title[Lane-Emden system in a half-space]{Liouville theorem for monotone solutions \\
of the Lane-Emden system in a half-space}
\author[Yimei Li and Philippe Souplet]{Yimei Li $^{(1)}$ and Philippe Souplet $^{(2)}$}

\thanks{$^{(1)}$School of Mathematics and Statistics, Beijing Jiaotong University, Beijing 100044, PR China.}

\thanks{$^{(2)}$Universit\'e Sorbonne Paris Nord, CNRS UMR 7539, LAGA, 93430 Villetaneuse, France.}

\thanks{Emails: liyimei@bjtu.edu.cn, souplet@math.univ-paris13.fr}

\date{}

\begin{document}
\begin{abstract}
The Lane-Emden system is an important model case of Hamitonian elliptic systems.
We here consider the zero-Dirichlet problem in a half-space for arbitrary powers $p,q\ge1$
and space dimension $n$. 
We prove that any nonnegative classical solution that is nondecreasing
in the normal direction is necessarily trivial.
No bound at infinity of any kind is assumed and
this extends the known scalar monotonicity theorem.

We also obtain the
same conclusion for the scalar equation $-\Delta u=f(u)$ when
$f\in C^1([0,\infty))$ is nonnegative, $f(0)=0$,
$f\not\equiv0$, and $f(s)/s$ is nondecreasing for $s>0$.
This relaxes the convexity and $C^2$ regularity 
requirements on $f$ in previous works.

To prove these results we introduce a new method, completely different from
the normal convexity approach from previous works,
which relies on a delicate combination of techniques
and could be applied to more general systems.
Our basic tool is a suitable average of the linear interpolant between the solution and its normal limit
on each finite strip.
We next use stability estimates, combined 
with an energy argument for the scalar equation, or a contraction
argument relying on Bessel potentials and the Schur test for the system. Increasing the strip height
then forces the solution to vanish.
\end{abstract}

\maketitle

\section{Introduction and main results}
\label{sec:introduction}

\subsection{Background}
The Lane-Emden system
\be{LERn}
 -\Delta u=v^p,\quad -\Delta v=u^q
\ee
is an important model case of Hamitonian elliptic system.
As such, it has received considerable interest in the literature (see references below)
and has also served as a benchmark for many methods in the analysis of nonlinear elliptic systems,
for instance variational or topological methods, Liouville type theorems, rescaling or moving planes arguments, etc.
In this paper, we are interested in the nonexistence problem for the half-space system
\begin{equation}\label{eq:le-problem}
 \begin{cases}
 -\Delta u=v^p&\text{in }\Rp,\\
 -\Delta v=u^q&\text{in }\Rp,\\
 u=v=0&\text{on }\partial\Rp,
 \end{cases}
\end{equation}
where $\Rp=\{x=(y,t)\in\R^{n-1}\times\R:t>0\}$, $n\ge2$,
and $p,q\ge1$. 
We also consider its scalar counterpart with a
general nonlinearity,
\begin{equation}\label{eq:scalar-problem}
 \begin{cases}
 -\Delta u=f(u)&\text{in }\Rp,\\
 u=0&\text{on }\partial\Rp.
 \end{cases}
\end{equation}
By a classical solution we mean 
$u\in C^2(\Rp)\cap C(\overline{\Rp})$ satisfying the equation and the
boundary condition pointwise, componentwise in the case of a system.
All solutions considered below are nonnegative.

To motivate our results, we first recall some background on the
corresponding problems in the whole space. 
For the Lane-Emden equation 
\[
 -\Delta u=u^p\quad\text{in }\R^n,
\]
the classical theorem of Gidas and Spruck~\cite{GS} gives nonexistence
of positive solutions when $1<p<p_S$, where
$p_S=(n+2)/(n-2)$ for $n\ge3$ and $p_S=\infty$ for $n\le2$.
For $n\ge3$ and $p\ge p_S$, positive radial bounded solutions exist.
Other proofs and related developments can be found in
\cite{bidaut1991nonlinear,chen1991classification,li2003liouville}.
For the  Lane-Emden system \eqref{LERn} 
in $\R^n$ the role of the Sobolev exponent is played by the Sobolev hyperbola
\[
 \frac1{p+1}+\frac1{q+1}=\frac{n-2}{n}
\]
(see \cite{clement1992positive,peletier1992existence}). Positive radial
bounded solutions exist on and above the critical hyperbola, that is,
when the left-hand side is at most $(n-2)/n$; see
\cite{serrin1998existence}. The Lane-Emden conjecture asks for
nonexistence in the opposite, subcritical range. 
Partial results were obtained in \cite{defigueiredo1994liouville,mitidieri1996nonexistence,serrin1996nonexistence,
reichel2000nonexistence,busca2002liouville,polacik2007singularity,souplet2009proof,LLW25}
covering in particular lower dimensions $n\le 4$.
Very recently, the preprint \cite{OpenAI26} proposed a proof of the full conjecture.

The situation in the half-space is different. Liouville theorems for
Dirichlet problems in a half-space are important both for their
intrinsic interest and for their applications to a priori estimates
through boundary rescaling arguments; see
\cite{gidas1981priori,quittner2019superlinear}. 
 They are in particular involved in the Berestycki-Caffarelli-Nirenberg conjecture
(cf.~\cite{BCN93,BCN98} and see, e.g.,~\cite{FV10, chen2014monotonicity,DF22}).
For the scalar case, Gidas and Spruck~\cite{gidas1981priori} proved nonexistence
for $1<p\le p_S$ by moving planes arguments. For bounded solutions,
Dancer~\cite{dancer1992some} obtained the improved range
$1<p<(n+1)/(n-3)$ for $n\ge4$, with no upper restriction when
$n\le3$. The argument uses $u_t>0$ in $\Rp$ and reduces the
problem to an equation in one dimension lower. 
By exploiting stability, Farina~\cite{farina2007classification}
obtained a further improved range for bounded solutions, covering in particular all dimensions $n\le 11$.
Then Chen, Lin and
Zou~\cite{chen2014monotonicity} introduced a normal convexity approach,
which allowed them to prove nonexistence of bounded
solutions for all $p>1$. More generally, their result applies to any nonlinearity $f\not\equiv 0$ such that
\be{hypCLZ}
\hbox{$f\ge 0$ convex,  with $f(0)=0$ and $f\in C^1([0,\infty))\cap C^2((0,\infty))$.}
\ee
For possibly unbounded scalar solutions, Dupaigne, Sirakov and
Souplet~\cite{DSS} proved that $u_t\ge0$ in $\Rp$ implies
nonexistence for every $p>1$, without a boundedness or growth
assumption on the solution.
Their proof
establishes normal convexity by a nonlinear maximum principle,
using a Moser iteration and stability estimates to control a weight
involving the normal derivative. 
 More generally, their result
applies to any nontrivial $f$ satisfying \eqref{hypCLZ},
and the normal monotonicity assumption $u_t\ge0$ is in particular 
satisfied whenever the solution is bounded on each finite strip
$\R^{n-1}\times(0,T)$, as a consequence of moving planes arguments. 
However it is unknown whether normal monotonicity
 is true for any classical solution (except in dimension $n=2$; see \cite{FS17}). 
Subsequently, for any $p>1$, Dupaigne, Farina and Petitt~\cite{dupaigne2023liouville},
enlarged the nonexistence class to 
classical solutions stable outside a compact set. However their proof crucially uses scale invariance
and their result does not cover more general nonlinearities.
Recently, Beuvin and Farina~\cite{BF26} obtained one dimensionality and nonexistence 
results for monotone solutions with nonconvex nonlinearities in dimensions $n\le 10$
(see Remark~\ref{rem:BF} for more details).

For the  Lane-Emden system \eqref{eq:le-problem}, early results on the half-space 
were obtained from Liouville theorems in $\mathbb R^n$.
In particular, Birindelli and Mitidieri~\cite{birindelli1998liouville} transferred
nonexistence of bounded whole-space solutions to the half-space in
one dimension higher. Reichel and Zou~\cite{reichel2000nonexistence}
proved nonexistence for possibly unbounded solutions in the range
$p,q<p_S$ by moving spheres arguments. Results based on stability
were obtained by Cowan~\cite{cowan2013liouville}; see also
\cite{CDG,DGH} for the whole-space stability theory. 
 By means of their normal convexity approach, Chen, Lin and
Zou~\cite{chen2014monotonicity} obtained nonexistence of bounded solutions
for every $p,q>1$
and also treated a suitable class of systems with convex nonlinearities.
More recently,
Li and Souplet~\cite{LS} proved nonexistence for every $p,q>1$ under
the assumption that the solution is bounded on each finite strip.
Their proof first establishes $u_t, v_t\ge 0$ in $\Rp$,
then derives logarithmic gradient estimates and normal convexity. The finite-strip boundedness assumption is used in the
gradient and maximum principle arguments, whereas no bound is
required as~$t\to\infty$.

This raises the question whether the scalar monotonicity theorem
of \cite{DSS} can be extended to systems without imposing any
boundedness or growth condition. For scalar equations,
a further question is whether the  convexity and $C^2$ regularity assumptions on $f$, 
which are crucial to the methods in \cite{chen2014monotonicity,DSS}, can be replaced by 
 weaker requirements.
The present
paper gives an affirmative answer to both questions, with
$u_t\ge0$ in the scalar case and $u_t,v_t\ge0$
in the system case, throughout $\Rp$, retained as hypotheses. The scalar result also
requires only $C^1$ regularity of $f$, possibly nonconvex, 
and both proofs avoid establishing convexity of $u$ 
in the normal direction.

\subsection{Main results}
Our first result concerns the system \eqref{eq:le-problem}. It shows that
$u_t,v_t\ge0$ forces $u=v=0$ without any bound on the solution.
\begin{theorem}\label{thm:system}
Let $n\ge2$ and $p,q\ge1$. If a nonnegative classical solution
$(u,v)\in C^2(\Rp)^2\cap C(\overline{\Rp})^2$ of
\eqref{eq:le-problem} satisfies $u_t\ge0$ and $v_t\ge0$ in $\Rp$,
then $u=v=0$.
\end{theorem}

The theorem extends
the scalar power result of \cite{DSS} to the Lane-Emden system and
replaces the finite-strip boundedness hypothesis of \cite{LS} by
 the weaker assumption
$u_t, v_t\ge0$ 
 (recall that the latter follows from finite-strip boundedness 
by moving planes arguments). 

\smallskip

For the scalar equation \eqref{eq:scalar-problem}, we obtain the same conclusion
under $u_t\ge0$, replacing convexity by the weaker condition $sf'(s)\ge f(s)$. 

\begin{theorem}\label{thm:scalar}
Let $n\ge2$ and let $f\in C^1([0,\infty))$ satisfy  $ f(0)=0$ and for any $s>0$,
\begin{equation}\label{eq:scalar-assumptions}
f(s)\ge0,\qquad f\not\equiv0,
 \qquad sf'(s)\ge f(s).
\end{equation}
If a nonnegative classical solution
$u\in C^2(\Rp)\cap C(\overline{\Rp})$ of
\eqref{eq:scalar-problem} satisfies $u_t\ge0$ in $\Rp$,
then $u=0$.
\end{theorem}

The last condition in \eqref{eq:scalar-assumptions} is equivalent
to the monotonicity of $f(s)/s$ on $(0,\infty)$. 
Every nonnegative convex $C^1$ function
with $f(0)=0$ satisfies this condition, but convexity is not
necessary. 
 For example, for all $p\ge 1$ and $k>k_0(p)>0$ sufficiently large, the function $f(s)=s^p+k(s^2-s\sin s)$
satisfies all the assumptions of  Theorem~\ref{thm:scalar}  but has a nonconvexity region around $s=3\pi/2$.

On the other hand, the conclusion of Theorem~\ref{thm:scalar} is in particular true if $u$ is assumed to be bounded on finite-strips instead of $u_t\ge0$, since this again implies the latter by moving planes arguments
(see \cite[Remark~3.1]{DSS} and the references therein).

\begin{remark} \label{rem:BF}
It is shown in \cite[Theorem~1.2]{BF26}  that $u_t\ge0$ in $\Rp$ implies one-dimensionality, i.e. $u(x)=u(t)$, 
whenever
\be{hypBF}
\hbox{$n\le 9$ and $0\le f\in C^1([0,\infty))$ satisfies } \displaystyle\liminf_{s\to\infty} \textstyle\frac{f(s)}{s}>0,
\ee
and that $u\equiv 0$ if, moreover, $f>0$ on $(0,\infty)$. 
For $n\le 9$, our Theorem~\ref{thm:scalar} is a consequence of~\cite[Theorem~1.2]{BF26} (by simple ODE analysis).
Theorem~\ref{thm:scalar} establishes this conclusion in every dimension.

Let us also mention that the same conclusions as in \cite[Theorem~1.2]{BF26} are obtained for $n=10$ in \cite[Theorem~1.4]{BF26} under the
additional assumption that $f$ has at most polynomial growth at infinity.
On the other hand, when $f\in C^1([0,\infty))$ is nonnegative nondecreasing and $\liminf_{s\to\infty}f(s)/s=\infty$,
 \cite[Theorems~1.6 and 1.9]{BF26} prove nonexistence of stable solutions for $n\le 4$, 
and of solutions stable outside a compact for $n\le 10$ if $f$ also has at most polynomial growth at infinity.
\end{remark}

\subsection{Plan of the proofs}
To prove these results we introduce a novel method, completely different from
the normal convexity approach from \cite{chen2014monotonicity,DSS,LS}. It relies on a delicate combination of techniques, 
which could be applied to more general systems.
In order not to obscure the novelty of the method with unessential technicalities, we have restricted ourselves to the 
important case of the Lane-Emden system and to the scalar case. 

Let us explain the main ideas. Write $d=n-1$ and
$x=(y,t)$. If $u\not\equiv0$ in the scalar case, or
$(u,v)\not\equiv(0,0)$ in the system case, the normal derivatives
are positive and provide stability inequalities 
(Propositions~\ref{prop:normal-stability-Scalar} and \ref{prop:normal-stability-system}).
In the scalar case, a uniform local energy
estimate follows from a first-eigenfunction test and the interior
gradient estimate of \cite{CFRS} (see Lemma~\ref{lem:scalar-energy}). 
For the system, the power scaling
and the stability inequality yield local bounds on $u^{q+1}+v^{p+1}$  (Lemma~\ref{lem:system-energy}).
These estimates pass to the limit and guarantee that the monotone normal limits
\[
 U(y)=\lim_{t\to\infty}u(y,t),\qquad
 V(y)=\lim_{t\to\infty}v(y,t)
\]
are finite a.e.~in $\R^d$ and satisfy polynomial integral
bounds (Propositions~\ref{prop:scalar-limit} and \ref{prop:system-limits}). 
These limits solve the corresponding stationary equation
or system in $\D'(\R^d)$. 

The existence of these limiting solutions does not by itself yield the desired conclusion.
On each finite strip, we use the positive part of the difference
between the original solution and its affine interpolant.
Let $\boldsymbol u, \boldsymbol U$ respectively denote the solution and its normal limit
(scalar or vector valued -- the same construction applies to scalar equations and systems).
On $\R^d\times(0,T)$, for each component, we define
\[
 W_{T,i}=\Big(u_i-\frac tT U_i\Big)_+,
 \qquad
 F_{T,i}(y)=\int_0^T\sin\Big(\frac{\pi t}{T}\Big)\,W_{T,i}(y,t)\dd t.
\]
The functions $\boldsymbol u, \boldsymbol U$ satisfy
\be{EqnsUu}
 -\Delta_y\bU=\bC\bU,\qquad -\Delta\bu\le\bC\bu, 
\ee
where the nonnegative coefficient matrix $\boldsymbol C(y)$ is given by
 $\bC=A=f(U)/U$ for the scalar equation, and by
$\bC=\begin{pmatrix}0&a\\b&0\end{pmatrix}$
with $(a,b)=(V^{p-1}(y),U^{q-1}(y))$ for the Lane-Emden system.
From \eqref{EqnsUu}, Kato's inequality and projection onto the first Dirichlet eigenfunction yield
\be{DeltaLambda}
 (-\Delta_y+\lambda_T)\bF_T\le\bC\bF_T,
 \qquad \lambda_T=\frac{\pi^2}{T^2}>0,
\ee
in the distributional sense (Proposition~\ref{lem:common-projection}).

The proofs that the projections vanish are then different. In the scalar
case, the positivity of the normal
limit and its local Sobolev regularity (a consequence of the above local energy estimate)  
yield the form boundedness property of the coefficient $A$: 
\be{DeltaLambda2}
 \int_{\R^d}Ah^2\dd y\le\int_{\R^d}|\nabla h|^2\dd y,\quad \quad h\in H^1(\R^d) 
\ee
(Lemma~\ref{prop:scalar-form}). 
Testing the projected inequality \eqref{DeltaLambda} with $F_T\eta^2$ and using 
\eqref{DeltaLambda2} 
 leaves
\[
 \lambda_T\int_{\R^d}F_T^2\eta^2\dd y
 \le\int_{\R^d}F_T^2|\nabla\eta|^2\dd y.
\]
A cutoff iteration and the polynomial integral bound imply $F_T=0$
(see Proposition~\ref{lem:scalar-extinction} and its proof). 

 For the system, this approach does not work 
because the form boundedness fails, due to the Hamiltonian coupling of the equations,
and a different idea is needed. 
Namely we establish a crucial resolvent contraction property
\be{MainResolv}
 \left\|\sqrt b\,(-\Delta_y+\lambda)^{-1}\sqrt a\right\|_{L^2(\R^d)\to L^2(\R^d)}
 \le\frac1{\sqrt{pq}}<1
\ee
(Proposition~\ref{prop:system-resolvent}).
The proof of \eqref{MainResolv} is delicate.
We first show a suitable, bounded cylinder version of \eqref{MainResolv}
(see Lemma~\ref{lem:finite-height}), 
by means of a two-weight Schur test (cf.~Proposition~\ref{lem:schur}).
To construct the required auxiliary functions for the test,
we use the positive functions $u_t$ and $v_t$, which satisfy the linearized system,
and consider an averaging of these functions in the normal direction against a sine function on a well-adjusted strip.
Here, the factors $p$ and $q$ in the equations for $u_t$ and $v_t$ produce the strict constant $1/\sqrt{pq}$.
The desired inequality \eqref{MainResolv} then follows by sending the height of the cylinder -- 
and next its radius -- to infinity, using monotone convergence of the corresponding Dirichlet resolvents
to the Bessel kernel (cf.~Proposition~\ref{rem:dirichlet-resolvent}).
Note that the above procedure requires
neither normalized auxiliary functions nor pointwise regularity
of $U$ and $V$. 
Then, although $\sqrt bF_{T,1}$ and
$\sqrt aF_{T,2}$ may not belong to $L^2(\R^d)$, their polynomial
integral growth allows us to multiply them by a decaying
exponential weight. An inequality between the Bessel kernels
at two nearby positive parameters (Lemma~\ref{lem:kernel-weight}) 
keeps the resulting operator norm
strictly smaller than one. The two projected inequalities  in \eqref{DeltaLambda} then
imply $F_{T,1}=F_{T,2}=0$
(Proposition~\ref{lem:system-extinction}). 

Finally, if these projections vanish, then
$u_i(y,t)\le tU_i(y)/T$ holds for $y\in\R^d\setminus\mathcal N_T$ and $0<t<T$, where $\mathcal N_T\subset\R^d$ and $|\mathcal N_T|=0$.
 Taking integer strip heights and letting $T\to\infty$
gives $\bu=0$ (Lemma~\ref{lem:common-finish}). 

\smallskip

The organization of the paper is as follows.
Section~\ref{sec:limits} establishes stability and local energy
estimates, and the normal limits. Section~\ref{sec:projection}
contains the common finite-strip projection and its final
consequence.  Theorems~\ref{thm:scalar} and~\ref{thm:system}
 are proved in Sections~\ref{sec:scalar-vanishing} and
\ref{sec:system-vanishing}, respectively.

\section{Local estimates  and normal limits}
\label{sec:limits}

We use $B_R(X)\subset\R^n$ for full-dimensional balls and
$B_R^d(y_0)=\{y\in\R^d:|y-y_0|<R\}$, 
and we write $B_R:=B_R(0)$ and $B_R^d:=B_R^d(0)$ for short. 
The Laplacian splits as $\Delta=\Delta_y+\partial_t^2$. Test functions in $C_c^1(\Rp)$
have support at positive distance from the flat boundary.
We write $H_c^1(\R^d)$ for the functions in $H^1(\R^d)$ with
compact support. The spaces $\D'(\Omega)$ and $\Sch'(\R^d)$
denote distributions and tempered distributions, respectively.
All constants denoted by $C$ are positive 
and may change from line
to line; their relevant dependence is stated in each result.

\subsection{Scalar local estimates and the normal limit}
We first strengthen $u_t\ge0$ to $u_t>0$ when $u\not\equiv0$.
 It is standard that positive derivative then yields the stability inequality.
We recall the short proof for convenience. 

\begin{proposition}\label{prop:normal-stability-Scalar}
Under the hypotheses of Theorem~\ref{thm:scalar}, if
$u\not\equiv0$, then $u>0$, $u_t>0$ in $\Rp$,
and $u$ is stable: For every
$\varphi\in C_c^1(\Rp)$, we have 
\begin{equation}\label{eq:scalar-stability}
 \int_{\Rp}f'(u)\varphi^2\dd x\le\int_{\Rp}|\nabla\varphi|^2\dd x. 
\end{equation}
\end{proposition}

\begin{proof}
Since $-\Delta u=f(u)\ge0$, the strong maximum principle gives
$u>0$ in $\Rp$. Set $\xi=u_t\ge0$. Differentiating the equation
and using local boundedness of the coefficient and interior $L^q$ elliptic regularity, we obtain a strong solution of
\[
 -\Delta\xi=f'(u)\xi\quad\text{in }\Rp.
\]
The strong maximum principle therefore shows that either
$\xi>0$ or $\xi\equiv0$. If $\xi\equiv0$, the zero boundary condition would imply $u=0$.
Hence $\xi>0$.

For $\varphi\in C_c^1(\Rp)$, 
using the identity
$$\xi^2\left|\nabla\left(\frac\varphi\xi\right)\right|^2
=\left|\nabla\varphi-\frac{\varphi\nabla\xi}{\xi}\right|^2
=|\nabla\varphi|^2-\nabla\xi\cdot\nabla\left(\frac{\varphi^2}{\xi}\right)
=|\nabla\varphi|^2+\frac{\Delta\xi}{\xi}\varphi^2
-\nabla\cdot\left(\frac{\varphi^2\nabla\xi}{\xi}\right)$$
we obtain
\begin{equation}\label{eq:ground-state-identity}
 \int_{\Rp}|\nabla\varphi|^2\dd x
 -\int_{\Rp}\frac{-\Delta\xi}{\xi}\varphi^2\dd x
 =\int_{\Rp}\xi^2\left|\nabla\left(\frac\varphi\xi\right)\right|^2\dd x
 \ge0.
\end{equation}
 Substitution of the equation for $\xi$ proves \eqref{eq:scalar-stability}.
\end{proof}

For the scalar case and $s>0$, set
\begin{equation}\label{eq:scalar-a}
 a_0(s)=\frac{f(s)}s.
\end{equation}
By \eqref{eq:scalar-assumptions}, the functions $a_0$ and $f$ are nondecreasing
and $f'(s)\ge a_0(s)\ge0$ for $s>0$. 
Moreover,
\begin{equation}\label{eq:scalar-secant}
 f(r)\le a_0(s)r,\quad s>0,\ 0\le r\le s.
\end{equation}
Also, choosing $s_*>0$ with $c_1:=f(s_*)/s_*>0$ and putting
$c_2=f(s_*)$,
we obtain
\begin{equation}\label{eq:linear-lower}
 f(s)\geq c_1s-c_2,\quad s\ge 0.
\end{equation}

 In the next lemma, after giving a basic a priori estimate, 
 we deduce local $H^1$~bounds independent of the solution,
as a consequence of the interior gradient estimate for stable solutions in
Cabr\'e--Figalli--Ros-Oton--Serra~\cite{CFRS}.
The latter is valid in every dimension and requires only a nonnegative,
locally Lipschitz nonlinearity. 

\begin{lemma}\label{lem:scalar-energy}
Assume \eqref{eq:scalar-assumptions}. There exist
$R_*=R_*(n,f)>1$ and $C=C(n,f)>0$ with the following properties.

\begin{itemize}[itemindent=-4mm]
\item[(i)] Every classical solution $u\ge 0$ of $-\Delta u=f(u)$ in $B_{R_*}(X)$ satisfies
\begin{equation}\label{eq:scalar-energy0}
 \int_{B_1(X)}f(u)\dd x \le C.
\end{equation}

\item[(ii)]  If, moreover, $u$ is stable, then
\begin{equation}\label{eq:scalar-energy}
 \int_{B_1(X)}\bigl(u^2+|\nabla u|^2\bigr)\dd x
 \le C.
\end{equation}
\end{itemize}
\end{lemma}

\begin{proof}
(i) After translation, we may take $X=0$.  
Denote by $\lambda_r$ the first Dirichlet eigenvalue of $B_r$, 
choose $r>8$ so large that $\lambda_r\le c_1/2$, and set $R_*=r+1$.
The closed ball $\overline{B_r}$ lies inside $B_{R_*}$, so $u$
is $C^2$ in a neighborhood of $\overline{B_r}$. Let $\varphi_1>0$ be the first
Dirichlet eigenfunction on $B_r$. Green's formula and 
$\partial_\nu\varphi_1\le 0=\varphi_1$ on $\partial B_r$ give
\[
 \int_{B_r}f(u)\varphi_1
\strut = -\int_{B_r}\varphi_1\Delta u
 =\lambda_r\int_{B_r}u\varphi_1
   +\int_{\partial B_r}u\,\partial_\nu\varphi_1
 \leq\lambda_r\int_{B_r}u\varphi_1.
\]
Together with \eqref{eq:linear-lower}, this yields
\be{eq:L1-new0}
c_1\int_{B_r}u\varphi_1
 \leq 2c_2\int_{B_r}\varphi_1.
\ee
Since $\varphi_1$ has a positive minimum on
$\overline{B_{r/2}}$, we conclude that
\begin{equation}\label{eq:L1-new}
 \int_{B_{r/2}}u\dd x\leq C.
\end{equation}
Choose $\phi\in C_c^\infty(B_{r/2})$, $0\leq\phi\leq1$,
with $\phi=1$ on $B_1$. Since $f(u)\geq0$, testing the equation
with $\phi$ and using \eqref{eq:L1-new} gives
$$  \int_{B_1}f(u)\dd x
\leq\int_{B_{r/2}}f(u)\phi\dd x
 =\int_{B_{r/2}}u(-\Delta\phi)\dd x
\strut\le  C\int_{B_{r/2}}u\dd x\leq C.
 $$ 
 
(ii) Applying \cite[Proposition~2.5]{CFRS}, after rescaling,
to the stable solution $u$, and using \eqref{eq:L1-new}, we obtain
\begin{equation}\label{eq:gradient-new}
 \int_{B_{r/4}}|\nabla u|^2\dd x\leq C.
\end{equation}
Set $\bar u=\fint_{B_{r/4}}u \dd x$. 
Applying Poincar\'e--Wirtinger's inequality and using
\eqref{eq:L1-new}, \eqref{eq:gradient-new}, we obtain
\begin{equation}\label{eq:L2-new}
\int_{B_{r/4}}u^2\le 
 2\int_{B_{r/4}}(u-\bar u)^2\dd x+2|B_{r/4}|\bar u^2\leq C \int_{B_{r/4}}|\nabla u|^2\dd x\
 +2|B_{r/4}|\bar u^2\leq C. 
\end{equation}
Combining \eqref{eq:gradient-new}-\eqref{eq:L2-new}, we obtain \eqref{eq:scalar-energy}.
\end{proof}

\begin{remark} \label{rem:linear}
The case $p=q=1$ of Theorem~\ref{thm:system}
 is an immediate consequence of the basic a priori estimate.
Indeed, for a nontrivial solution, since $w=u+v$ satisfies $-\Delta w=w$ and $w>0$, 
inequality \eqref{eq:L1-new0} with $c_1=1$ and $c_2=0$
in \eqref{eq:linear-lower} leads to a contradiction.
\end{remark}

The preceding local bounds remain uniform when the solution is translated
upward. We use them to show that the normal limit is finite a.e.~and satisfies $-\Delta_yU=f(U)$ in $\D'(\R^d)$.

\begin{proposition}\label{prop:scalar-limit}
Let $u\not\equiv0$ satisfy the hypotheses of Theorem~\ref{thm:scalar}. The limit $U(y)=\lim_{t\to\infty}u(y,t)$
is finite and positive a.e.~in $\R^d$, belongs to
$H^1_{\loc}(\R^d)$, and satisfies
\begin{equation}\label{eq:scalar-limit-equation}
 -\Delta_yU=f(U)\quad\text{in }\D'(\R^d).
\end{equation}
For every $R\ge1$,
\begin{equation}\label{eq:scalar-limit-energy}
 \int_{B_R^d}(U^2+f(U))\dd y\le C(1+R)^d.
\end{equation}
\end{proposition}

\begin{proof}
Fix $R\geq1$ and $t>0$. For every $y\in B_R^d$ and
$\tau\in[t+R_*+1,t+R_*+2]$, the inequality $u_t\geq0$ gives
\[
 u(y,t)^2\leq u(y,\tau)^2.
\]
Thus
\begin{equation}\label{eq:horizontal-slice}
 \begin{split}
 \int_{B_R^d}u(y,t)^2\dd y
 &=\int_{t+R_*+1}^{t+R_*+2}\int_{B_R^d}u(y,t)^2\dd y\dd\tau\\
 &\leq\int_{B_R^d\times[t+R_*+1,t+R_*+2]}u(y,\tau)^2\dd y\dd\tau.
 \end{split}
\end{equation}
Choose balls $B_1(X_j)$ covering this cylinder such that
\[
 \begin{gathered}
 B_R^d\times[t+R_*+1,t+R_*+2]
 \subset\displaystyle\bigcup_{j=1}^N B_1(X_j),
 \qquad N\leq C_n(1+R)^d,\\
 (X_j)_n\in[t+R_*+1,t+R_*+2].
 \end{gathered}
\]
Since $(X_j)_n\geq t+R_*+1>R_*$, each
$B_{R_*}(X_j)\Subset\Rp$.
Lemma~\ref{lem:scalar-energy} therefore gives
\begin{equation}\label{eq:cylinder-bound}
 \begin{split}
 \int_{B_R^d\times[t+R_*+1,t+R_*+2]}u^2\dd x
 &\leq\sum_{j=1}^N\int_{B_1(X_j)}u^2\dd x\leq C N\leq C(1+R)^d,
 \end{split}
\end{equation}
where $C$ is independent of $R$ and $t$.
Combining \eqref{eq:horizontal-slice}--\eqref{eq:cylinder-bound}
and letting $t\to\infty$, monotone convergence gives
\[
 \int_{B_R^d}U(y)^2\dd y
 =\lim_{t\to\infty}\int_{B_R^d}u(y,t)^2\dd y
 \leq C(1+R)^d.
\]
Thus $U$ is finite a.e.~in each $B_R^d$, hence a.e.~in
$\R^d$ by a countable exhaustion. Since $f$ is nonnegative,
nondecreasing and continuous, the same covering argument, now
using the bound in \eqref{eq:scalar-energy0}, gives
\[
 \int_{B_R^d}f(u(y,t))\dd y\leq C(1+R)^d.
\]
Letting $t\to\infty$ and applying monotone convergence yields
the bound for $f(U)$ in \eqref{eq:scalar-limit-energy}.

Let $D\subset\R^d$ be a ball and put $K=D\times(0,1)$.
For integers $j$ sufficiently large, the translates
$u_j(y,t)=u(y,j+t)$ are uniformly bounded in $H^1(K)$ by
Lemma~\ref{lem:scalar-energy} and a finite covering of $K$.
As $j\to\infty$, we have $u_j\to U$ a.e.~in $K$. Since $0\le u_j\le U$ and
$U\in L^2(D)$, dominated convergence gives $u_j\to U$
in $L^2(K)$ as $j\to\infty$. Weak compactness therefore identifies the
$t$-independent function $U(y)$ as an element of $H^1(K)$.
Averaging its weak tangential derivatives over $(0,1)$ shows
that $U\in H^1(D)$, hence $U\in H^1_{\loc}(\R^d)$.

Similarly, $0\le f(u_j)\le f(U)\in L^1(D)$ and continuity
of $f$ give $f(u_j)\to f(U)$ in $L^1(K)$ as $j\to\infty$. For
$\varphi\in C_c^\infty(D)$ and $\rho\in C_c^\infty((0,1))$
with $\int_0^1\rho=1$, test the  equation with
$\varphi(y)\rho(t)$ and let $j\to\infty$. Since
$\int_0^1\rho''=0$, this yields
\[
 \int_D U(-\Delta_y\varphi)\dd y
 =\int_D f(U)\varphi\dd y.
\]
Thus \eqref{eq:scalar-limit-equation} holds. Finally,
$U(y)\ge u(y,1)>0$ a.e.~in $\R^d$ by
Proposition~\ref{prop:normal-stability-Scalar}.
\end{proof}

\subsection{System local estimates and the normal limits}

In the system estimates assume $p,q\ge1$, $pq>1$, and set
\begin{equation}\label{eq:scaling}
 \alpha=\frac{2(p+1)}{pq-1},\qquad
 \beta=\frac{2(q+1)}{pq-1},\qquad
 s=\alpha(q+1)=\beta(p+1)=\alpha+\beta+2.
\end{equation}

For the system, positivity of both derivatives gives a stability inequality
with the coefficient $\sqrt{pq}\,u^{(q-1)/2}v^{(p-1)/2}$.
A similar result was given in \cite[Lemma~3]{cowan2013liouville} for stable solutions of \eqref{LERn} in $\R^n$.
We give the short proof for convenience. 

\begin{proposition}\label{prop:normal-stability-system}
	Under the hypotheses of Theorem~\ref{thm:system}, if
$(u,v)\not\equiv(0,0)$, then
\[
 u,v,u_t,v_t>0\quad\text{in }\Rp.
\]
For every $\varphi\in C_c^1(\Rp)$, we also have
	\begin{equation}\label{eq:system-stability}
		\sqrt{pq}\int_{\Rp}u^{(q-1)/2}v^{(p-1)/2}\varphi^2\dd x
		\le\int_{\Rp}|\nabla\varphi|^2\dd x.
	\end{equation}
\end{proposition}

\begin{proof}
If one component vanished identically, the equations would force
the other to vanish as well. Hence the strong maximum principle
gives $u,v>0$, and interior regularity implies that $u,v$ are
smooth. Set $\xi=u_t\ge0$, $\eta=v_t\ge0$. Differentiation
gives 
\begin{equation}\label{eq:jacobi}
 -\Delta\xi=pv^{p-1}\eta,\qquad
 -\Delta\eta=qu^{q-1}\xi.
\end{equation}
If either derivative vanished identically, both would vanish by
these equations, and the zero boundary values would imply
$u=v=0$. Since $\xi,\eta$ are nonnegative and superharmonic,
the strong maximum principle therefore gives $\xi,\eta>0$.

Let $\varphi\in C_c^1(\Rp)$. 
The  identity \eqref{eq:ground-state-identity},
applied to each derivative, yields
\[
 \int_{\Rp}pv^{p-1}\frac\eta\xi\varphi^2\dd x
 \le\int_{\Rp}|\nabla\varphi|^2\dd x,
 \qquad
 \int_{\Rp}qu^{q-1}\frac\xi\eta\varphi^2\dd x
 \le\int_{\Rp}|\nabla\varphi|^2\dd x.
\]
Adding the inequalities and using
\[
\sqrt{pq}\,u^{(q-1)/2}v^{(p-1)/2}\le 
\frac12\Big(pv^{p-1}\frac\eta\xi+qu^{q-1}\frac\xi\eta\Big)
\]
gives \eqref{eq:system-stability}. 
\end{proof}

The next estimate converts the stability inequality into bounds for
$u^{q+1}+v^{p+1}$. 
A similar result is given in \cite[Section~3]{DGH}, assuming a global comparison of the components.
We here provide a modified proof which does not require this property. 
These bounds will control the 
 normal limits and their nonlinearities
(only the case $\theta=1$ will actually be used).

\begin{lemma}\label{lem:system-energy}
Let $p,q\ge1$ with $pq>1$. 
Let $(u,v)$ be a positive classical solution of \eqref{LERn} 
in $B_{4R}(X)$, satisfying
\eqref{eq:system-stability} for every $\varphi\in C_c^1(B_{4R}(X))$,
with both integrals taken over $B_{4R}(X)$. 
For $1\le \theta\le n/(n-2)$ 
if $n\ge3$, or any finite $\theta\ge 1$ if
$n=2$, there is $C=C(n,p,q,\theta)>0$ such that
\begin{equation}\label{eq:system-higher-energy}
 \int_{B_R(X)}(u^{q+1}+v^{p+1})^\theta\dd x
 \le CR^{n-s\theta}.
\end{equation}
\end{lemma}

\begin{proof}
Constants in this proof depend only on $n,p,q,\theta$ and the chosen cutoffs.
Also, in view of H\"older's inequality, we may assume $\theta>1$ without loss of generality. 

\smallskip
\textbf{Step 1.} \textit{Scaling.}
Set $\widetilde u(z)=R^\alpha u(X+Rz)$ and
$\widetilde v(z)=R^\beta v(X+Rz)$. By \eqref{eq:scaling}, the identities
$\alpha+2=p\beta$, $\beta+2=q\alpha$ preserve the equations,
and $\alpha(q-1)+\beta(p-1)=4$ preserves stability.
It is therefore enough to prove the estimate for $R=1$, $X=0$.

\smallskip
\textbf{Step 2.} \textit{Integral bounds for $u^q$ and $v^p$.}
Let $\varphi>0$ be the first Dirichlet eigenfunction on
$B_{7/2}$, normalized by $\int_{B_{7/2}}\varphi=1$, and let
$\Lambda$ be its eigenvalue. 
Integrating by parts and using
$\partial_\nu\varphi\le 0$ on $\partial B_{7/2}$ gives
\[
 \int_{B_{7/2}}v^p\varphi\le\Lambda\int_{B_{7/2}}u\varphi,
 \qquad
 \int_{B_{7/2}}u^q\varphi\le\Lambda\int_{B_{7/2}}v\varphi.
\]
Writing $m_u=\int_{B_{7/2}}u\varphi$ and
$m_v=\int_{B_{7/2}}v\varphi$, Jensen's inequality yields
$m_v^p\le\Lambda m_u$ and $m_u^q\le\Lambda m_v$. Hence
$m_u^{pq-1}\le\Lambda^{p+1}$ and 
$m_v^{pq-1}\le\Lambda^{q+1}$.
Since $pq>1$ and $\varphi$ has a positive lower bound on $B_3$,
the preceding estimates imply
\begin{equation}\label{eq:nonlinear-integral-bound}
 \int_{B_3}(u^q+v^p)\dd x\le C.
\end{equation}

\smallskip
\textbf{Step 3.} \textit{The cross term and Sobolev improvement.}
Put
\[
 S=u^{q+1}+v^{p+1},\qquad D=u^{(q-1)/2}v^{(p-1)/2},
 \qquad\gamma=\frac{q+1}{2},\quad\delta=\frac{p+1}{2},
\]
and
\[
 A=\sqrt{pq}\frac q{\gamma^2},\qquad
 B=\sqrt{pq}\frac p{\delta^2}.
\]
Since $2p\ge p+1$ and $2q\ge q+1$, with at least one strict inequality owing to $pq>1$, we have 
\begin{equation}\label{eq:energy-strict-bound}
 \sqrt{AB}=\frac{4pq}{(p+1)(q+1)}>1.
\end{equation}
For $0\le\eta\in C_c^2(B_2)$, define
\[
 J=\int_{B_2}u^qv^p\eta^2\dd x,\qquad
 I_1=\int_{B_2}D u^{q+1}\eta^2\dd x,\qquad
 I_2=\int_{B_2}D v^{p+1}\eta^2\dd x,
\]
\[
 E_u=\int_{B_2}u^{q+1}(|\nabla\eta|^2+|\Delta(\eta^2)|)\dd x,
 \quad
 E_v=\int_{B_2}v^{p+1}(|\nabla\eta|^2+|\Delta(\eta^2)|)\dd x,
 \quad E=E_u+E_v.
\]
Testing $-\Delta u=v^p$ with $u^q\eta^2$ and integrating by
parts gives
\[
 q\int_{B_2}u^{q-1}|\nabla u|^2\eta^2\dd x
 =J+\frac1{q+1}\int_{B_2}u^{q+1}\Delta(\eta^2)\dd x.
\]
Expanding $\nabla(u^\gamma\eta)$ and integrating its mixed
term by parts therefore yields
\be{SobolevEstim}
 \int_{B_2}|\nabla(u^\gamma\eta)|^2\dd x
 =\frac{\gamma^2}{q}J
 +\int_{B_2}u^{q+1}|\nabla\eta|^2\dd x
 +\frac{1-\gamma}{2q}\int_{B_2}u^{q+1}\Delta(\eta^2)\dd x
\ee
and an analogous calculation on $-\Delta v=u^q$ gives
\be{SobolevEstim2}
 \int_{B_2}|\nabla(v^\delta\eta)|^2\dd x
 =\frac{\delta^2}p J
 +\int_{B_2}v^{p+1}|\nabla\eta|^2\dd x
 +\frac{1-\delta}{2p}\int_{B_2}v^{p+1}\Delta(\eta^2)\dd x.
\ee
The stability inequality \eqref{eq:system-stability} applied with $\varphi=u^\gamma\eta$
and $\varphi=v^\delta\eta\in C_c^1(B_2)$ 
then yields
\[
 AI_1\le J+CE_u,\qquad BI_2\le J+CE_v.
\]
Observing that $(Du^{q+1}Dv^{p+1})^{1/2}=u^qv^p$,
we get $J\le\sqrt{I_1I_2}$ by the Cauchy--Schwarz inequality.
It follows that
\[
 \sqrt{AB}\,J\le\sqrt{(J+CE_u)(J+CE_v)}\le J+CE
\]
and \eqref{eq:energy-strict-bound} then implies 
\begin{equation}\label{eq:cross-estimate}
 J\le CE.
\end{equation}
Going back to \eqref{SobolevEstim}-\eqref{SobolevEstim2}, 
we obtain
\[
 \int_{B_2}\bigl(|\nabla(u^\gamma\eta)|^2+
                    |\nabla(v^\delta\eta)|^2\bigr)\dd x
 \le C(J+E)\le CE.
\]
Since $S\eta^2=(u^\gamma\eta)^2+(v^\delta\eta)^2$, the
Sobolev inequality in $H_0^1(B_2)$
yields
\begin{equation}\label{eq:poisson-improvement}
 \|S\eta^2\|_{L^\theta(B_2)}
 \le\|u^\gamma\eta\|_{L^{2\theta}(B_2)}^2
     +\|v^\delta\eta\|_{L^{2\theta}(B_2)}^2
 \le CE.
\end{equation}

\smallskip
\textbf{Step 4.} \textit{Interpolation and conclusion.} 
Let $\kappa=\min\{q/(q+1),p/(p+1)\}\in(0,1)$.
Since $\kappa(q+1)\le q$, $\kappa(p+1)\le p$, and
$s^r\le1+s^k$ for $0<r\le k$, the inequality
$(a+b)^\kappa\le a^\kappa+b^\kappa$ and H\"older's inequality,
followed by \eqref{eq:nonlinear-integral-bound}, give
\begin{equation}\label{eq:low-moment}
 \int_{B_2}S^\kappa\dd x
 \le\int_{B_2}\bigl(u^{\kappa(q+1)}+v^{\kappa(p+1)}\bigr)\dd x
 \le C.
\end{equation}
In view of \eqref{eq:low-moment}, choose $\chi\in C_c^\infty(B_2)$ with $0\le\chi\le1$ and
$\chi=1$ on $B_1$, and set $\eta=\chi^m$ for an integer $m$.
Writing
\[
 J_1=\int_{B_2}S^\theta\chi^{2m\theta}\dd x,\qquad
 J_2=\int_{B_2}S\chi^{2m-2}\dd x,
\]
we have $J_1^{1/\theta}\le C_mJ_2$ by
\eqref{eq:poisson-improvement}. Set
\[
 \ell=\frac{\theta-1}{\theta-\kappa}\in(0,1),\qquad
 1=\kappa\ell+\theta(1-\ell).
\]
Since $\theta(1-\ell)<1$, we may choose $m$ sufficiently large so that
$2m[1-\theta(1-\ell)]-2\ge0$. H\"older's inequality then gives
\[
 J_2
 \le J_1^{1-\ell}
 \left(\int_{B_2}S^\kappa
 \chi^{\{2m[1-\theta(1-\ell)]-2\}/\ell}\dd x\right)^\ell
 \le CJ_1^{1-\ell}\le CJ_2^{\theta(1-\ell)},
\]
hence $J_2\le C$ and then $J_1\le C$.
This implies 
$\int_{B_1}S^\theta\le C$ and completes the proof. 
\end{proof}

We can now pass the local system bounds to horizontal slices and then let
the height $t\to\infty$. The next proposition identifies the resulting limiting equations.

\begin{proposition}\label{prop:system-limits}
Let $(u,v)$ satisfy the hypotheses of Theorem~\ref{thm:system},
with $(u,v)\not\equiv(0,0)$ and $pq>1$. The limits
$U(y)=\lim_{t\to\infty}u(y,t)$ and
$V(y)=\lim_{t\to\infty}v(y,t)$ are finite and positive almost
everywhere in $\R^d$.
For every $R>0$,
\begin{equation}\label{eq:system-limit-energy}
 \int_{B_R^d}(U^{q+1}+V^{p+1})\dd y\le CR^{d-s}.
\end{equation}
In particular, $U,V,U^q,V^p$ are locally integrable.
They solve
\begin{equation}\label{eq:system-limit-equations}
 -\Delta_yU=V^p,\qquad -\Delta_yV=U^q
 \quad\text{in }\D'(\R^d).
\end{equation}
\end{proposition}

\begin{proof}
\textbf{Step 1.} \textit{Integral bounds and finiteness.}
Fix $R>0$, $t>0$, and put $X=(0,t+9R)$. The ball $B_{8R}(X)$
lies inside $\Rp$. Proposition~\ref{prop:normal-stability-system}
provides its stability inequality. Applying
Lemma~\ref{lem:system-energy} with $\theta=1$ gives
\[
 \int_{B_{2R}(X)}(u^{q+1}+v^{p+1})\dd x\le CR^{n-s}.
\]
The cylinder $B_R^d\times(t+8R,t+10R)$ is contained in
$B_{2R}(X)$. Since $u_t>0$ and $v_t>0$ in $\Rp$, we have
\[
 2R\int_{B_R^d}(u(y,t)^{q+1}+v(y,t)^{p+1})\dd y
 \le\int_{B_R^d\times(t+8R,t+10R)}(u^{q+1}+v^{p+1})\dd x
 \le CR^{n-s}.
\]
Dividing by $2R$ and letting $t\to\infty$, monotone convergence
yields \eqref{eq:system-limit-energy}. In particular, the limits
are finite a.e.~in each $B_R^d$, hence a.e.~in
$\R^d$ by a countable exhaustion. The functions are positive since
$U(y)\ge u(y,1)>0$ and $V(y)\ge v(y,1)>0$ on the set where they are finite, outside a null set.
H\"older's inequality on bounded balls gives the local
integrability of $U,V,U^q,V^p$ from
\eqref{eq:system-limit-energy}, also when $p=1$ or $q=1$.

\smallskip
\textbf{Step 2.} \textit{The limiting equations.}
Let $T_j\to\infty$ and set
$u_j(y,t)=u(y,T_j+t)$, $v_j(y,t)=v(y,T_j+t)$. On any fixed
compact cylinder $K\times[-L,L]$, these translates are defined
for $j$ large. As $j\to\infty$, we have
\[
 0\le u_j\le U,\quad 0\le v_j\le V,
 \qquad u_j\to U,\quad v_j\to V\quad\text{a.e.~in }K\times[-L,L].
\]
Since $U,V,U^q,V^p\in L^1(K)$, letting $j\to\infty$ and applying
dominated convergence, we obtain
\[
 u_j\to U,\quad v_j\to V,\quad
 u_j^q\to U^q,\quad v_j^p\to V^p
 \quad\text{in }L^1(K\times[-L,L]).
\]
Take $\varphi\in C_c^\infty(\R^d)$ and
$\rho\in C_c^\infty(\R)$ with $\int_\R\rho=1$. Testing the
first equation with $\varphi(y)\rho(t)$ gives
\[
 \int u_j[-\rho\Delta_y\varphi-\varphi\rho'']\dd y\dd t
 =\int v_j^p\varphi\rho\dd y\dd t.
\]
Letting $j\to\infty$ and using $\int_\R\rho''=0$, we obtain
\[
 \int_{\R^d}U(-\Delta_y\varphi)\dd y
 =\int_{\R^d}V^p\varphi\dd y.
\]
The second equation follows in the same way. This proves
\eqref{eq:system-limit-equations}.
\end{proof}

\section{A common finite-strip reduction}\label{sec:projection}

A main tool of the proofs of Theorems~\ref{thm:system} and \ref{thm:scalar}
 is the following coercive elliptic inequality for suitable normal averages 
of the difference between the original solution and its affine interpolant
  on a finite strip.
  It will be eventually applied to a solution and its normal limit.
  The arguments in this section are componentwise and do not use
the particular scalar or Lane-Emden nonlinearity. 
All vector inequalities below are understood componentwise.

\begin{proposition}\label{lem:common-projection}
Let $m\ge1$.  Let $\bu=(u_1,\ldots,u_m)$, $\bU=(U_1,\ldots,U_m)$
with $u_i\in C^2(\Rp)\cap C^1_{\loc}(\overline{\Rp})$ and $U_i\in L^1_{loc}(\R^d)$ satisfy 
$$\hbox{$\bu(y,0)=0\le\bu(y,t)\le\bU(y)$ a.e.~in $\Rp$,\ \ for all $t>0$.}$$
Let $\bC(y)=(c_{ij}(y))$ be a matrix on $\R^d$ with a.e.~finite nonnegative, measurable 
entries such that $\bC\bU\in L^1_{\loc}(\R^d)^m$, and assume
\begin{equation}\label{eq:common-hypotheses}
 -\Delta_y\bU=\bC\bU\quad\text{in }\D'(\R^d)^m,
 \qquad -\Delta\bu\le\bC\bu\quad\text{in }\D'(\Rp)^m.
\end{equation}
For $T>0$, set $F_{T,i}(y)=\displaystyle\int_0^T\psi_T(t)W_{T,i}(y,t)\dd t$, where
\begin{equation}\label{eq:common-positive-parts}
 w_{T,i}=u_i-\frac tT U_i,\qquad W_{T,i}=(w_{T,i})_+,\\
 \qquad \psi_T(t)=\sin\frac{\pi t}{T},\qquad
 \lambda_T=\frac{\pi^2}{T^2}.
\end{equation}
Then $\bF_T=(F_{T,1},\ldots,F_{T,m})^{\mathsf T}\in L^\infty_{\loc}(\R^d)^m$ satisfies
\begin{equation}\label{eq:common-projection}
 (-\Delta_y+\lambda_T)\bF_T\le\bC\bF_T
 \quad\text{in }\D'(\R^d)^m.
\end{equation}
Moreover, for $1\le i\le m$, we have
\begin{equation}\label{eq:common-domination}
 0\le F_{T,i}\le\frac{2T}{\pi}U_i
 \quad\text{a.e.~in }\R^d
\end{equation}
and, if $U_i\in H^1_{\loc}(\R^d)$ then
$F_{T,i}\in H^1_{\loc}(\R^d)$.
\end{proposition}

In view of the proof we state the following distributional version of Kato's inequality.
Since the assumptions are weaker than those usually encountered in the literature,
we give a proof in appendix for convenience.

\begin{lemma}\label{lem:kato}
Let $\Omega\subset\R^m$ be open and $w,g\in L^1_{\loc}(\Omega)$.
If $\Delta w\ge g$ in $\D'(\Omega)$, then
\begin{equation}\label{eq:kato}
 \Delta w_+\ge\ind_{\{w>0\}}g\quad\text{in }\D'(\Omega).
\end{equation}
\end{lemma}

\begin{proof}[Proof of Proposition~\ref{lem:common-projection}]
Using Fubini's theorem and continuity in $t$, we may choose a measurable
set $E\subset\R^d$ with $|\R^d\setminus E|=0$. For $y\in E$,
$t\ge0$ and $1\le i\le m$, we have
\[
 0\le u_i(y,t)\le U_i(y)<\infty.
\]
We use this set throughout the proof.

\smallskip
\textbf{Step 1.} \textit{Equations for the positive parts.} 
Fix $T>0$ and write $\theta=t/T$, $w_i=u_i-\theta U_i$,
$W_i=(w_i)_+$. Since $\partial_t^2(\theta U_i)=0$, the assumptions
\eqref{eq:common-hypotheses} imply
\[
 -\Delta\bw=-\Delta\bu+\theta\Delta_y\bU
 \le\bC\bu-\theta\bC\bU=\bC\bw
 \quad\text{in }\D'(S_T)^m,
\]
where $\bw=(w_1,\ldots,w_m)^{\mathsf T}$ and
$\Delta=\Delta_y+\partial_t^2$.
All terms on the right are locally integrable, because
$|w_j|\le2U_j$ and $\bC\bU\in L^1_{\loc}$. Applying
Lemma~\ref{lem:kato} to each component gives
\begin{equation}\label{eq:common-kato}
 -\Delta W_i
 \le\ind_{\{w_i>0\}}\sum_{j=1}^m c_{ij}w_j
 \le\sum_{j=1}^m c_{ij}W_j
 \quad\text{in }\D'(S_T).
\end{equation}
The last inequality uses $c_{ij}\ge0$.

\smallskip
\textbf{Step 2.} \textit{Bounds at the strip endpoints.}
Let $K\subset\R^d$ be compact. Since $T$ is fixed, local
boundary regularity gives a constant $C_{K,T}$ such that
$W_i(y,t)\le u_i(y,t)\le C_{K,T}t$ on $K\times(0,T)$.
For $y\in K\cap E$, the inequality $U_i(y)\ge u_i(y,T)$ gives
\[
 \begin{split}
 W_i(y,t)
 &\le\Big(u_i(y,t)-\frac tT u_i(y,T)\Big)_+
 \le |u_i(y,t)-u_i(y,T)|
       +\Big(1-\frac tT\Big)u_i(y,T)
 \le C_{K,T}(T-t).
 \end{split}
\]
Consequently, for $y\in K\cap E$ and $0<t<T$, we have
\begin{equation}\label{eq:common-endpoint-bound}
 0\le W_i(y,t)\le C_{K,T}\min\{t,T-t\}.
\end{equation}
The inequality $W_i\le u_i$ also gives
$F_{T,i}\in L^\infty_{\loc}(\R^d)$.

\smallskip
\textbf{Step 3.} \textit{Projection and removal of the endpoint cutoffs.}
For $0<\varepsilon<T/4$, choose
$\eta_\varepsilon\in C_c^\infty(0,T)$ with
$0\le\eta_\varepsilon\le1$, $\eta_\varepsilon=1$ on
$[2\varepsilon,T-2\varepsilon]$, and
\[
 |\eta_\varepsilon'|\le C\varepsilon^{-1},\qquad
 |\eta_\varepsilon''|\le C\varepsilon^{-2}.
\]
Its derivatives are supported in
$(0,2\varepsilon)\cup(T-2\varepsilon,T)$. For
$0\le\varphi\in C_c^\infty(\R^d)$,
test \eqref{eq:common-kato} with
$\varphi(y)\psi_T(t)\eta_\varepsilon(t)$. Since
$-\psi_T''=\lambda_T\psi_T$, we obtain
\[
 \begin{split}
 &\int_{S_T}W_i(-\Delta_y\varphi+\lambda_T\varphi)
                   \psi_T\eta_\varepsilon\dd y\dd t-\int_{S_T}W_i\varphi
            (2\psi_T'\eta_\varepsilon'+\psi_T\eta_\varepsilon'')\dd y\dd t\\
 &\le\int_{S_T}\sum_jc_{ij}W_j\varphi\psi_T\eta_\varepsilon\dd y\dd t.
 \end{split}
\]
Set $K=\operatorname{supp}\varphi$,
$I_\varepsilon=(0,2\varepsilon)\cup(T-2\varepsilon,T)$ and
$r(t)=\min\{t,T-t\}$. For $y\in K\cap E$ and $0<t<T$, we have
\[
 W_i(y,t)\le C_{K,T}r(t),\qquad
 0\le\psi_T(t)\le C_Tr(t),\qquad |\psi_T'(t)|\le C_T.
\]
Since $\eta_\varepsilon',\eta_\varepsilon''$ vanish outside
$I_\varepsilon$, we obtain as $\varepsilon\downarrow0$,
\[
 \begin{aligned}
\left|\int_{S_T}W_i\varphi
       (2\psi_T'\eta_\varepsilon'+\psi_T\eta_\varepsilon'')\dd y\dd t\right|
 &\le\int_K\int_{I_\varepsilon}W_i(y,t)\varphi(y)
       \bigl(2|\psi_T'||\eta_\varepsilon'|
                 +\psi_T|\eta_\varepsilon''|\bigr)\dd t\dd y\\
 &\le C\int_{I_\varepsilon}r(t)
                 (\varepsilon^{-1}+r(t)\varepsilon^{-2})\dd t\\
 &=2C\int_0^{2\varepsilon}t
                 (\varepsilon^{-1}+t\varepsilon^{-2})\dd t
       \le C\varepsilon\to 0.
 \end{aligned}
\]
Here $C$ is independent of $\varepsilon$.
As $\varepsilon\downarrow0$, the other integrals converge by dominated convergence:
on the right, $\sum_jc_{ij}W_j\le(\bC\bU)_i$. Letting
$\varepsilon\downarrow0$ and using Fubini's theorem gives
\[
 \int_{\R^d}F_{T,i}(-\Delta_y\varphi+\lambda_T\varphi)\dd y
 \le\int_{\R^d}\sum_jc_{ij}F_{T,j}\varphi\dd y.
\]
This is \eqref{eq:common-projection}. The domination follows from
$W_i\le U_i$ and $\int_0^T\psi_T\dd t=2T/\pi$.

\smallskip
\textbf{Step 4.} \textit{Local Sobolev regularity when $U_i\in H^1_{\loc}$.}
For any ball $D=B_r^d(y_0)\subset\R^d$, we have
$w_i\in H^1(D\times(0,T))$, with
\[
 \nabla_yw_i=\nabla_yu_i-\frac tT\nabla_yU_i,
 \qquad \partial_tw_i=(u_i)_t-\frac1T U_i.
\]
The positive-part chain rule gives $W_i\in H^1$ on this cylinder.
Consequently,
\[
 \partial_{y_\ell}F_{T,i}
 =\int_0^T\psi_T(t)\partial_{y_\ell}W_i(\cdot,t)\dd t
 \quad\text{in }\D'(D),
\]
and Cauchy--Schwarz bounds this derivative in $L^2(D)$. This
proves $F_{T,i}\in H^1_{\loc}(\R^d)$ and completes the proof.
\end{proof}

In the next sections, by two different methods applicable to the scalar problem and to the system,
respectively, we will show that inequality \eqref{eq:common-projection} forces 
the projections $\bF_T$ on each strip to vanish. 
By the following simple lemma, this will then imply the triviality of the solution.

\begin{lemma}\label{lem:common-finish}
Under the assumptions of Proposition~\ref{lem:common-projection},
suppose $\bF_k=0$ a.e.~in $\R^d$ for every 
integer $k\ge 1$. Then $\bu=0$ in $\Rp$.
\end{lemma}

\begin{proof}
For each integer $k\ge1$, choose a null set $\mathcal N_k$
outside which all $U_i$ are finite, $F_{k,i}=0$, and
$u_i(y,t)\le U_i(y)$ for every $t>0$. Such a choice is possible
by Fubini's theorem and the continuity of $t\mapsto u_i(y,t)$.
For $y\notin\mathcal N_k$, we have 
\[
 0=F_{k,i}(y)
 =\int_0^k\sin\frac{\pi t}{k}
        \left(u_i(y,t)-\frac tk U_i(y)\right)_+\dd t.
\]
The integrand is nonnegative and the sine is strictly positive
on $(0,k)$. Hence
$(u_i(y,t)-tU_i(y)/k)_+=0$ almost
everywhere in $(0,k)$, and then everywhere by continuity. Thus,
for $0<t<k$, we have
\[
 0\le u_i(y,t)\le\frac tk U_i(y).
\]
Let $\mathcal N=\bigcup_{k\ge1}\mathcal N_k$. For every
$y\notin\mathcal N$ and $t>0$, letting $k\to\infty$ through
integers $k>t$ gives $u_i(y,t)=0$. The complement of
$\mathcal N$ is dense in $\R^d$. Continuity of $u_i(\cdot,t)$
then implies $u_i=0$ everywhere, for each component.
\end{proof}

\section{Proof of Theorem~\ref{thm:scalar}}

\label{sec:scalar-vanishing}

\subsection{Vanishing of the scalar projection by energy inequality}

The following proposition concerns the coercive inequality 
\begin{equation}\label{eq:scalar-subsolution}
 (-\Delta_y+\lambda)F\le AF
\end{equation} 
in the scalar case
(without direct reference to the original problem \eqref{eq:scalar-problem}).
It guarantees the vanishing of solutions under two assumptions: 
the form boundedness of~$A$ and a polynomial integral bound on $F$.
It is proved by an energy inequality and a cutoff iteration.

\begin{proposition}
\label{lem:scalar-extinction}
Let $\lambda>0$, $A\in L^1_{\loc}(\R^d)$ and
$F\in H^1_{\loc}(\R^d)$, with $A, F\ge 0$, satisfy
\eqref{eq:scalar-subsolution} in $\D'(\R^d)$.
Let $C>0$, $N\ge0$ and assume
\begin{equation}\label{eq:scalar-projection-growth}
 \int_{B_R^d}F^2\dd y\le C(1+R)^N,\quad R\ge1,
\end{equation}
and
\begin{equation}\label{eq:scalar-form}
 \int_{\R^d}Ah^2\dd y\le\int_{\R^d}|\nabla h|^2\dd y,\quad h\in H^1(\R^d).
\end{equation}
Then $F=0$ a.e.~in $\R^d$.
\end{proposition}

\begin{proof}
We first justify testing the distributional inequality with
$F\eta^2$, then use a cutoff iteration.

\smallskip

\textbf{Step 1.} \textit{Local integrability and the energy inequality.}
We claim that
\begin{equation}\label{eq:scalar-massive-cutoff}
 \lambda\int_{\R^d}F^2\eta^2\dd y
 \le\int_{\R^d}F^2|\nabla\eta|^2\dd y,
 \quad\hbox{for all $\eta\in C_c^\infty(\R^d)$ with $\eta\ge 0$.}
\end{equation}

For every compact $K\subset\R^d$, choose
\[
 \chi_K\in C_c^\infty(\R^d),\qquad
 0\le\chi_K\le1,\qquad \chi_K|_K=1.
\]
Since $F\chi_K\in H^1(\R^d)$, \eqref{eq:scalar-form} gives
\[
 \int_K AF^2\dd y
 \le\int_{\R^d}A(F\chi_K)^2\dd y
 \le\int_{\R^d}|\nabla(F\chi_K)|^2\dd y<\infty.
\]
By Cauchy--Schwarz inequality, we have
\[
 \int_K AF\dd y
 \le\left(\int_K A\dd y\right)^{1/2}
      \left(\int_K AF^2\dd y\right)^{1/2}<\infty.
\]
Thus $AF^2,AF\in L^1_{\loc}(\R^d)$.
Fix $0\le\eta\in C_c^\infty(\R^d)$ and choose $R>0$ with
$\operatorname{supp}\eta\subset B_R^d$. Set
\[
 h=F\eta^2,\qquad K_\eta=\overline{B_{R+1}^d}.
\]
Choose $0\le\rho\in C_c^\infty(B_1^d)$ with $\int_{\R^d}\rho\dd y=1$.
For integers $k\ge1$, set
\[
 \rho_k(y)=k^d\rho(ky),\qquad h_k=h*\rho_k.
\]
As $k\to\infty$, we have
 \[
 h_k\in C_c^\infty(\R^d),\qquad h_k\ge0,\qquad
 \operatorname{supp}h_k\cup\operatorname{supp}h\subset K_\eta,
 \qquad h_k\to h\quad\text{in }H^1(\R^d).
\]
Choose $\chi\in C_c^\infty(\R^d)$ with $0\le\chi\le1$ and
$\chi|_{K_\eta}=1$. Applying Cauchy--Schwarz inequality and \eqref{eq:scalar-form}, we obtain,
as $k\to\infty$,
\be{convh1}
 \begin{aligned}
 \left|\int_{\R^d}AF(h_k-h)\dd y\right|
 &\le\left(\int_{\R^d}A(F\chi)^2\dd y\right)^{1/2}
       \left(\int_{\R^d}A(h_k-h)^2\dd y\right)^{1/2}\\
 &\le\|\nabla(F\chi)\|_2\,\|\nabla(h_k-h)\|_2\longrightarrow0.
 \end{aligned}
\ee
Also, as $k\to\infty$, we have
\be{convh2}
 \left|\int_{\R^d}\nabla F\cdot\nabla(h_k-h)\dd y\right|
 \le\|\nabla F\|_{L^2(K_\eta)} \|\nabla(h_k-h)\|_2\longrightarrow0,
\ee
 and
\be{convh3}
 \left|\int_{\R^d}F(h_k-h)\dd y\right| \le\|F\|_{L^2(K_\eta)}\,\|h_k-h\|_2\longrightarrow0.
\ee

Now, testing \eqref{eq:scalar-subsolution} with $h_k\ge0$, we have 
\[
 \int_{\R^d}\nabla F\cdot\nabla h_k\dd y
 +\lambda\int_{\R^d}Fh_k\dd y
 \le\int_{\R^d}AFh_k\dd y.
\]
Using \eqref{convh1}-\eqref{convh3} we may let $k\to\infty$ to get 
$$ \int_{\R^d}\nabla F\cdot\nabla(F\eta^2)\dd y
 +\lambda\int_{\R^d}F^2\eta^2\dd y
 \le\int_{\R^d}A(F\eta)^2\dd y,
$$
and \eqref{eq:scalar-form} then implies
$$ \int_{\R^d}\nabla F\cdot\nabla(F\eta^2)\dd y
 +\lambda\int_{\R^d}F^2\eta^2\dd y \le\int_{\R^d}|\nabla(F\eta)|^2\dd y.$$
In view of the identity 
$$|\nabla(F\eta)|^2=\nabla F\cdot\nabla(F\eta^2)+F^2|\nabla\eta|^2,$$
inequality \eqref{eq:scalar-massive-cutoff} follows. 

\smallskip
\textbf{Step 2.} \textit{Iteration and extinction.}
To apply \eqref{eq:scalar-massive-cutoff}, we choose $\eta_R\in C_c^\infty(B_{2R}^d)$, $0\le\eta_R\le1$,
with $\eta_R=1$ on $B_R^d$ and $|\nabla\eta_R|\le c/R$.
Writing $E(R)=\int_{B_R^d}F^2\dd y$, we obtain, for $R\ge1$,
\[
 E(R)\le\frac{c^2}{\lambda R^2}E(2R).
\]
Iteration gives, for every integer $k\ge1$,
\[
 E(R)\le\left(\frac{c^2}{\lambda R^2}\right)^k
            2^{-k(k-1)}E(2^kR).
\]
For fixed $R$, using \eqref{eq:scalar-projection-growth}, we have
$E(2^kR)\le C_R2^{kN}$. Hence, as $k\to\infty$, we obtain
\[
 0\le E(R)\le C_R
       \left(\frac{c^2}{\lambda R^2}\right)^k
       2^{-k(k-1)+kN}\longrightarrow0.
\]
Thus $E(R)=0$ for every $R\ge1$, and $F=0$ a.e.~in $\R^d$.
\end{proof}

\subsection{Form boundedness of the coefficient $A$}

Let $u\not\equiv0$ satisfy the hypotheses of Theorem~\ref{thm:scalar}
and let $U(y)$ be the scalar normal limit given by Proposition~\ref{prop:scalar-limit}.
By modifying $U$ on a null set, we may assume without loss of generality 
that 
\be{Ufinite}
0<U(y)<\infty\quad\hbox{for all $y\in\R^d$}
\ee
 and we set
\be{eq:scalar-A}
A(y)= \dfrac{f(U(y))}{U(y)}.
\ee
The next lemma guarantees the form boundedness of $A$, 
needed to apply Proposition~\ref{lem:scalar-extinction}.
It is proved by testing equation \eqref{eq:scalar-limit-equation} with $h^2/U$.

\begin{lemma}
\label{prop:scalar-form}
The function $A$ in \eqref{eq:scalar-A} satisfies property \eqref{eq:scalar-form}.
\end{lemma}

\begin{proof}
\par
Let $h\in C_c^\infty(\R^d)$. Choose a bounded open set
$D\Subset\R^d$ containing its support. Positivity and continuity
of $u(\cdot,1)$ give
\[
 U\ge u(\cdot,1)\ge c_D>0\quad\text{a.e.~on }D.
\]
Since $U\in H^1(D)$, the function $\zeta=h^2/U$, extended
by zero outside $D$, belongs to $H_c^1(\R^d)\cap L^\infty$,
and
\[
 \nabla\zeta=\frac{2h}{U}\nabla h-\frac{h^2}{U^2}\nabla U.
\]
Nonnegative smooth mollifications $\zeta_k$ have support in a
fixed compact set, are uniformly bounded, and converge as $k\to\infty$
to $\zeta$ in $H^1$ and, along a subsequence, a.e.~in $\R^d$.
Since $f(U)\in L^1_{\loc}(\R^d)$, letting $k\to\infty$ and applying
dominated convergence, we obtain
\[
 \int_{\R^d}f(U)\zeta_k\dd y
 \longrightarrow\int_{\R^d}f(U)\zeta\dd y.
\]
Since $\nabla U\in L^2_{\loc}(\R^d)$ and
$\zeta_k\to\zeta$ in $H^1(\R^d)$ as $k\to\infty$, we have
\[
 \int_{\R^d}\nabla U\cdot\nabla\zeta_k\dd y
 \longrightarrow\int_{\R^d}\nabla U\cdot\nabla\zeta\dd y.
\]
 We may thus test the weak equation \eqref{eq:scalar-limit-equation} for $U$ with $\zeta$.
We obtain
\[
 \begin{split}
 \int_{\R^d}Ah^2\dd y
 &=\int_{\R^d}f(U)\frac{h^2}{U}\dd y
 =\int_{\R^d}\nabla U\cdot\nabla\left(\frac{h^2}{U}\right)\dd y\\
 &=\int_{\R^d}|\nabla h|^2\dd y
 -\int_{\R^d}\left|\nabla h-\frac hU\nabla U\right|^2\dd y
 \le\int_{\R^d}|\nabla h|^2\dd y.
 \end{split}
\]
For $h\in H_c^1(\R^d)$, choose $h_\ell\in C_c^\infty(\R^d)$
with supports in a fixed compact set and $h_\ell\to h$ in $H^1(\R^d)$
and a.e.~in $\R^d$ as $\ell\to\infty$. Letting $\ell\to\infty$
and applying Fatou's lemma, we obtain the inequality for $h$. Finally, for $h\in H^1(\R^d)$, choose
$\chi_R\in C_c^\infty(B_{2R}^d)$ with $0\le\chi_R\le1$,
$\chi_R=1$ on $B_R^d$, and $|\nabla\chi_R|\le C/R$.
As $R\to\infty$, we have $\chi_Rh\to h$ in $H^1$ and a.e.~in $\R^d$.
Letting $R\to\infty$ and applying Fatou's lemma, we obtain \eqref{eq:scalar-form}.
\end{proof}

\subsection{Completion of proof of Theorem~\ref{thm:scalar}}

Assume for contradiction that $u\not\equiv0$ satisfies the hypotheses of Theorem~\ref{thm:scalar}.
Let $U(y)$, satisfying \eqref{Ufinite}, be the scalar normal limit given by Proposition~\ref{prop:scalar-limit}
and let $A$ be given by \eqref{eq:scalar-A}.
For every compact set $K\subset\R^d$, positivity and continuity
give $\min_Ku(\cdot,1)>0$. Since $U\ge u(\cdot,1)$, we get
\[
 A\le\frac{f(U)}{\min_Ku(\cdot,1)}\quad\text{a.e.~on }K.
\]
Proposition~\ref{prop:scalar-limit} therefore gives
$A\in L^1_{\loc}(\R^d)$. Moreover,
$-\Delta_yU=f(U)=AU$ in $\D'(\R^d)$. For all $t>0$ we have
$0\le u(y,t)\le U(y)$ a.e.~$y\in\R^d$, so \eqref{eq:scalar-secant} implies
$f(u)\le Au$. Thus $(m,\bu,\bU,\bC)=(1,u,U,A)$ satisfies all the assumptions of Proposition~\ref{lem:common-projection},
including $U\in H^1_{\loc}(\R^d)$.
It follows that, for each $T>0$, $F_T\in H^1_{\loc}(\R^d)\cap L^\infty_{\loc}(\R^d)$ 
and $\lambda_T:=\frac{\pi^2}{T^2}$ satisfy
$$
 (-\Delta_y+\lambda_T)F_T\le AF_T\quad\text{in }\D'(\R^d),
 \qquad 0\le F_T\le\frac{2T}{\pi}U\quad\text{a.e.~in }\R^d.
$$
In view of \eqref{eq:scalar-limit-energy}, the latter inequality implies \eqref{eq:scalar-projection-growth} with $N=d$.
Moreover,
\eqref{eq:scalar-form} is satisfied owing to Lemma~\ref{prop:scalar-form}.
We may thus apply Proposition~\ref{lem:scalar-extinction} to deduce, for every $T>0$, that
$F_T=0$ a.e.~in $\R^d$, hence 
$u\equiv0$ in $\Rp$ by Lemma~\ref{lem:common-finish}.
This contradicts $u\not\equiv0$ and concludes the proof of Theorem~\ref{thm:scalar}.
\qed

\section{Proof of Theorem~\ref{thm:system}}
\label{sec:system-vanishing}

By means of Proposition~\ref{lem:common-projection}, 
for suitable weights $\bar a, \bar b\ge 0$, the problem will be reduced to showing 
the vanishing of projections $F_1, F_2$ which satisfy a system of the form
\be{SystF1F2}
\left\{
 \begin{aligned}
 (-\Delta_y+\lambda)F_1&\le\bar a F_2\\
(-\Delta_y+\lambda)F_2&\le\bar b F_1.
 \end{aligned}
 \right.
 \ee
 
\subsection{Resolvents, polynomial integral growth, Bessel kernels and Schur test} 

To tackle system \eqref{SystF1F2}, we need various preliminary properties of the resolvent $R_\lambda=(-\Delta_y+\lambda)^{-1}$
 and of the corresponding Bessel kernels, that we collect in this subsection. 
We first fix the resolvent notation and the class of functions whose
absolute integrals grow at most polynomially. This class will allow us to use kernels without global integrability.

\begin{definition}\label{def:resolvent}
For $\lambda>0$ let $R_\lambda=(-\Delta_y+\lambda)^{-1}$ be
the Fourier multiplier $(|\xi|^2+\lambda)^{-1}$ on
$\Sch'(\R^d)$. Its positive Bessel kernel is
\begin{equation}\label{eq:bessel-kernel}
 G_\lambda(z)=\int_0^\infty e^{-\lambda r}(4\pi r)^{-d/2}
                      e^{-|z|^2/(4r)}\dd r.
\end{equation}
Also let
\begin{equation}\label{eq:Pclass}
 \Pclass=\left\{f\in L^1_{\loc}(\R^d):
 \int_{B_R^d}|f|\dd y\le C(1+R)^N\ (R\ge1)
 \text{ for some }C>0,\ N\ge0\right\}.
\end{equation}
\end{definition}

The kernel has mass $1/\lambda$ and satisfies
$(-\Delta_y+\lambda)G_\lambda=\delta_0$. On $L^2$, the
resolvent is self-adjoint. Products with nonnegative
multiplication coefficients will always be interpreted by
their positive integral kernels.

The following four propositions are essentially known, so we only give proofs in Appendix~A for the convenience of readers.
The first one shows that the resolvent preserves the class $\Pclass$
defined in \eqref{eq:Pclass}. 
In particular, the convolution remains well defined locally even when $f$ is not globally integrable.

\begin{proposition}\label{lem:resolvent-P}
We have
\begin{equation}\label{eq:resolvent-P}
 \Pclass\subset\Sch'(\R^d),\qquad R_\lambda(\Pclass)\subset\Pclass.
\end{equation}
For $f\in\Pclass$,  $R_\lambda f:=G_\lambda*f$ is locally
integrable, and satisfies
$(-\Delta_y+\lambda)R_\lambda f=f$ in $\mathcal{D}'(\R^d)$. 
\end{proposition}

To pass from estimates on balls to estimates on $\R^d$, we also need
the convergence of the Dirichlet kernels as the radius tends to infinity.

\begin{proposition}\label{rem:dirichlet-resolvent}
For a ball $D\subset\R^d$, let $R_{\lambda,D}$ denote the Dirichlet
resolvent of $-\Delta_y+\lambda$, with positive kernel $G_{\lambda,D}$.
Extend $G_{\lambda,B_R^d}$ by zero outside $B_R^d$ in both variables.
For fixed $x,z\in\R^d$ with $x\ne z$, we then have
\begin{equation}\label{eq:dirichlet-exhaustion}
 G_{\lambda,B_R^d}(x,z)\nearrow G_\lambda(x-z)
 \quad\text{as }R\to\infty.
\end{equation}
\end{proposition}

 The third tool is a maximum principle, of Phragmen-Lindel\"of type,
for the operator $-\Delta_y+\lambda$ in the class $\Pclass$.

\begin{proposition}\label{lem:resolvent-upper-bound}
Let $\lambda>0$ and $F,f\in\Pclass$. If
$(-\Delta_y+\lambda)F\le f$ in $\D'(\R^d)$, then
\begin{equation}\label{eq:resolvent-upper-bound}
 F\le R_\lambda f\quad\text{a.e.~in }\R^d.
\end{equation}
\end{proposition}

We next recall the two-weight Schur test. It turns two inequalities
for a positive kernel into an operator norm bound on $L^2(D)$.

\begin{proposition}\label{lem:schur}
	Let $D\subset\mathbb R^d$ be bounded and let $K$ be a nonnegative measurable kernel on $D\times D$. Suppose that there are measurable functions $h,j>0$ a.e.~on $D$ and constants $A,B>0$ such that
	\begin{equation}\label{eq:schur-hypotheses}
		\int_D K(x,z)h(z)\,dz\leq A\,j(x),
		\quad
		\int_D K(x,z)j(x)\,dx\leq B\,h(z),
		\quad\hbox{for a.e.~$x,z\in D$.}
	\end{equation}
	Then the integral operator $(Tf)(x)=\int_DK(x,z)f(z)\,dz$ is bounded on $L^2(D)$ and
	\begin{equation}\label{eq:schur-conclusion}
		\|T\|_{2\to2}:=\sup_{\substack{f\in L^2\\ f\neq0}}\frac{\|Tf\|_{L^2}}{\|f\|_{L^2}}\leq\sqrt{AB}.
	\end{equation}
\end{proposition}

Finally we record the following simple lemma, concerning exponential weighting on the Bessel kernel
in terms of shift of the parameter $\lambda$.
It gives a quantified estimate by a factor that tends to one as $\delta\downarrow0$, so as to preserve strict norm bounds.

\begin{lemma}
\label{lem:kernel-weight}
Let $0<\delta<1$ and $\varepsilon>0$ satisfy
$\varepsilon^2<\delta\lambda$. Set
\[
 C_\delta=(1-\delta)^{1-d/2},\qquad
 \lambda'=(1-\delta)(\lambda-\varepsilon^2/\delta)>0.
\]
Then
\begin{equation}\label{eq:kernel-weight}
 e^{\varepsilon|z|}G_\lambda(z)\le C_\delta G_{\lambda'}(z).
\end{equation}
\end{lemma}

\begin{proof}
Young's inequality gives, for $r=|z|$ and $t>0$,
\[
 \varepsilon r-\lambda t-\frac{r^2}{4t}
 \le-\left(\lambda-\frac{\varepsilon^2}{\delta}\right)t
     -\frac{(1-\delta)r^2}{4t}.
\]
Insert this estimate into
\eqref{eq:bessel-kernel} and use $t=(1-\delta)s$.
We then obtain
\[
 \begin{split}
 e^{\varepsilon|z|}G_\lambda(z)
 &\le\int_0^\infty(4\pi t)^{-d/2}
 e^{-(\lambda-\varepsilon^2/\delta)t-(1-\delta)|z|^2/(4t)}\dd t\\
 &=(1-\delta)^{1-d/2}
 \int_0^\infty(4\pi s)^{-d/2}e^{-\lambda's-|z|^2/(4s)}\dd s
 =C_\delta G_{\lambda'}(z).
 \end{split}
\]
The  inequality $\varepsilon^2<\delta\lambda$ ensures $\lambda'>0$.
\end{proof}

\subsection{Vanishing of the projections under contraction assumption}

In this section, for given measurable weights $\bar a,\bar b\ge 0$, we consider
a family of operators which appear naturally in connection with system \eqref{SystF1F2}.
For $\mu>0$, the operator $T_\mu:=\sqrt{\bar b} R_\mu\sqrt{\bar a}$ is defined on $L^1_{loc}(\R^d)$ by its nonnegative kernel, namely 
\be{defTmu}
(T_\mu f)(y):=\int_{\R^d} \sqrt{\bar b(y)}G_\mu(y-z) \sqrt{\bar a(z)}f(z)\,dz\in [0,\infty],
\quad 0\le f\in L^1_{loc}(\R^d).
\ee 
Whenever it is well defined, its extension to $L^2(\R^d)$ (by separation of the positive and negative parts of $f$)
will be still denoted by  $T_\mu$. 
The following proposition guarantees the vanishing of solutions of \eqref{SystF1F2} under 
 a contraction assumption on $T_\mu$ and suitable integral polynomial growth conditions.

\begin{proposition}\label{lem:system-extinction}
Let $\lambda>0$, $\theta\in(0,1)$, let $\bar a,\bar b\ge 0$ be measurable functions on $\R^d$ and assume that,
for every~$\mu\in(0,\lambda)$,
\be{hypTmuL2}
\|T_\mu\|_{L^2(\R^d)\to L^2(\R^d)}\le \theta.
\ee
Suppose that $0\le F_1,F_2\in L^1_{\loc}(\R^d)$ 
satisfy \eqref{SystF1F2} in $\D'(\R^d)$, along with
\be{newcondaF}
F_1, \, F_2, \, \bar bF_1, \, \bar aF_2,\, \bar bF_1^2,\, \bar aF_2^2\, \in\Pclass.
\ee
Then $F_1=F_2=0$ a.e.~in $\R^d$.
\end{proposition}

\begin{proof}
\smallskip
\textbf{Step 1.} \textit{Resolvent inequalities and weighted integrability.}
We set $T^*_\mu=\sqrt{\bar a}R_\mu\sqrt{\bar b}$,
defined similarly as in \eqref{defTmu}, and note that
\be{TmustarL2}
\|T^*_\mu\|_{L^2(\R^d)\to L^2(\R^d)}=\|T_\mu\|_{L^2(\R^d)\to L^2(\R^d)}
\ee
by $L^2$ duality. Namely, this follows from the identity
\[
 \int [g\,T_\mu f](y)\dd y
 =\int g(y)\sqrt{\bar b(y)}G_\mu(y-z)
           \sqrt{\bar a(z)}f(z)\dd z\dd y
 =\int [f\,T^*_\mu g](z)\dd z, 
\]
for all nonnegative $f,g\in L^2$, owing to Tonelli's theorem.

Next, by assumption \eqref{newcondaF},
we may apply Proposition~\ref{lem:resolvent-upper-bound} to deduce
\[
 F_1\le R_\lambda(\bar aF_2),\qquad
 F_2\le R_\lambda(\bar bF_1).
\]
Putting $X=\sqrt{\bar b}F_1$, $Y=\sqrt{\bar a}F_2$,
this rewrites as 
\begin{equation}\label{eq:system-resolvent-pair}
 X\le T_\lambda Y,\qquad Y\le T_\lambda^*X.
\end{equation}
Moreover, for $\varepsilon>0$ and $w_\varepsilon(y)=e^{-\varepsilon|y|}$, we have
\be{wespL2}
w_\varepsilon X, w_\varepsilon Y\in L^2(\R^d),
\ee
owing to $X^2,Y^2\in\Pclass$ (by \eqref{newcondaF}) and 
\[
 \int_{\R^d}w_\varepsilon^2X^2\dd y
 \le\sum_{j\ge0}e^{-2\varepsilon j}
                      \int_{B_{j+1}^d}X^2\dd y
 \le C\sum_{j\ge0}e^{-2\varepsilon j}(j+2)^N<\infty.
\]

\smallskip
\textbf{Step 2.} \textit{Strict contraction after weighting.}
Since $C_\delta=(1-\delta)^{1-d/2}\to1$ as $\delta\downarrow0$
we may choose $\delta\in(0,1)$ such that
$\kappa=C_\delta\theta <1$. Then choose
$0<\varepsilon<\sqrt{\delta\lambda}$, and let
$\lambda'=(1-\delta)(\lambda-\varepsilon^2/\delta)>0$.
The triangle inequality gives
$w_\varepsilon(x)/w_\varepsilon(z)\le e^{\varepsilon|x-z|}$.
Applying Lemma~\ref{lem:kernel-weight}, 
we obtain
\[
 \frac{w_\varepsilon(x)}{w_\varepsilon(z)}
   \sqrt{\bar b(x)}G_\lambda(x-z)\sqrt{\bar a(z)}
\le C_\delta\sqrt{\bar b(x)}G_{\lambda'}(x-z)\sqrt{\bar a(z)}.
\]
The same bound holds with $\bar a,\bar b$ interchanged.
Applying these kernel inequalities to
\eqref{eq:system-resolvent-pair}~gives
\[
 w_\varepsilon X\le C_\delta T_{\lambda'}(w_\varepsilon Y),
 \qquad
 w_\varepsilon Y\le C_\delta T_{\lambda'}^*(w_\varepsilon X).
\]
Taking norms and using \eqref{hypTmuL2}, \eqref{TmustarL2}, \eqref{wespL2}, it follows that 
\[
 \|w_\varepsilon X\|_2\le\kappa\|w_\varepsilon Y\|_2,
 \qquad
 \|w_\varepsilon Y\|_2\le\kappa\|w_\varepsilon X\|_2.
\]
Consequently, $\|w_\varepsilon X\|_2\le
\kappa^2\|w_\varepsilon X\|_2$, so both weighted norms vanish.
 Since $w_\varepsilon>0$, we have $X=Y=0$ a.e. in $\R^d$. Thus,
$\bar bF_1=\sqrt{\bar b}\,X=0$ and $\bar aF_2=\sqrt{\bar a}\,Y=0$ a.e.~in $\R^d$.
Substituting these identities into \eqref{SystF1F2} yields
\[
 (-\Delta_y+\lambda)F_i\le0
 \quad\text{in }\D'(\R^d),\qquad i=1,2.
\]
As $F_i\in\Pclass$, Proposition~\ref{lem:resolvent-upper-bound}
gives $F_i\le R_\lambda0=0$. Since $F_i\ge0$,
we conclude that $F_1=F_2=0$ a.e.~in $\R^d$.
\end{proof}

\subsection{Contraction property for Lane-Emden system}

Fix $p,q\ge 1$, let $(u,v)\not\equiv(0,0)$ satisfy the hypotheses of Theorem~\ref{thm:system},
and let $U, V$ be the normal limits in Proposition~\ref{prop:system-limits}. We set
$$a=V^{p-1},\ \ b=U^{q-1}\qquad\hbox{and}\qquad a_\tau=v^{p-1}(\cdot,\tau),\ \ b_\tau=u^{q-1}(\cdot,\tau)
\quad\hbox{for $\tau>0$.}$$
The next result guarantees the contraction property of the operator $T_\lambda=\sqrt b R_\lambda\sqrt a$, needed to apply Proposition~\ref{lem:system-extinction}.

\begin{proposition}\label{prop:system-resolvent}
The kernel
\begin{equation}\label{eq:weighted-resolvent-kernel}
 K_\lambda(x,z)=\sqrt{b(x)}G_\lambda(x-z)\sqrt{a(z)}
\end{equation}
defines a bounded operator
$T_\lambda=\sqrt b R_\lambda\sqrt a$ on $L^2(\R^d)$.
For every $\lambda>0$, we have
\begin{equation}\label{eq:system-resolvent-bound}
 \|T_\lambda\|_{2\to2}\le(pq)^{-1/2}.
\end{equation}
\end{proposition}

The proof will rely on a bounded domain approximation, reducing to the following key lemma, 
where $R_{\lambda,D}$ is the Dirichlet resolvent (cf.~Proposition~\ref{rem:dirichlet-resolvent}).

\begin{lemma}\label{lem:finite-height}
For any $\lambda,\tau>0$ and ball $D\subset\R^d$, there holds
\begin{equation}\label{eq:finite-height-bound}
 \left\|\sqrt{b_\tau}\,R_{\lambda,D}\sqrt{a_\tau}\right\|_{2\to2}
 \le(pq)^{-1/2}.
\end{equation}
\end{lemma}

The proof of Lemma~\ref{lem:finite-height}
uses $u_t,v_t>0$ to construct suitable auxiliary functions to apply the Schur test. 
Their sine averages give a resolvent bound independent of the starting height.

\begin{proof}
Fix $\lambda,\tau>0$ and a ball $D\subset\R^d$.

\smallskip
\textbf{Step 1.} \textit{Normal averaging.}
Applying Proposition~\ref{prop:normal-stability-system}, we have
$\xi=u_t>0$ and $\eta=v_t>0$, and these functions solve \eqref{eq:jacobi}.
Put $H=\pi/\sqrt\lambda$, $\psi(t)=\sin(\sqrt\lambda(t-\tau))$
on $[\tau,\tau+H]$, and define
\[
 \Phi(y)=\int_\tau^{\tau+H}\xi(y,t)\psi(t)\dd t,\qquad
 \Psi(y)=\int_\tau^{\tau+H}\eta(y,t)\psi(t)\dd t.
\]
Multiplying the first equation  of \eqref{eq:jacobi} 
by $\psi(t)$ and integrating over $t\in[\tau,\tau+H]$ gives,
using $\psi(\tau)=\psi(\tau+H)=0$,
$\psi'(\tau)=\sqrt\lambda$, $\psi'(\tau+H)=-\sqrt\lambda$,
and $\psi''=-\lambda\psi$,
\[
 \begin{aligned}
 \int_\tau^{\tau+H}\partial_t^2\xi(y,t)\psi(t)\dd t
 &=-\big[\xi(y,t)\psi'(t)\big]_\tau^{\tau+H}
   +\int_\tau^{\tau+H}\xi(y,t)\psi''(t)\dd t\\
 &=\sqrt\lambda\,[\xi(y,\tau)+\xi(y,\tau+H)]-\lambda\Phi(y).
 \end{aligned}
\]
Therefore,
\be{IneqDeltaPhi}
 (-\Delta_y+\lambda)\Phi(y)
 =p\int_\tau^{\tau+H}v(y,t)^{p-1}\eta(y,t)\psi(t)\dd t
  +\sqrt\lambda\,[\xi(y,\tau)+\xi(y,\tau+H)].
\ee
Since $v_t>0$ and $p\ge1$, for $y\in D$ and $\tau\le t\le\tau+H$,
we have
$v(y,t)^{p-1}\ge v(y,\tau)^{p-1}=a_\tau(y).$
For $y\in D$, we also have
$\sqrt\lambda\,[\xi(y,\tau)+\xi(y,\tau+H)]\ge0.$
Inequality \eqref{IneqDeltaPhi} and a similar calculation for the second equation in \eqref{eq:jacobi} thus give,
in the classical sense: 
\begin{equation}\label{ineqtaupsi}
 (-\Delta_y+\lambda)\Phi\ge p a_\tau\Psi,\qquad
 (-\Delta_y+\lambda)\Psi\ge q b_\tau\Phi.
\end{equation}
Consequently, $Z:=\Phi-pR_{\lambda,D}(a_\tau\Psi)$ 
satisfies $-\Delta_yZ+\lambda Z\ge0$ in $D$ with $Z=\Phi>0$ on $\partial D$,
hence $Z\ge 0$ in $D$ by the maximum principle.
Arguing similarly on $\Psi$ we thus obtain 
\begin{equation}\label{ClaimD}
 pR_{\lambda,D}(a_\tau\Psi)\le\Phi,\qquad
 qR_{\lambda,D}(b_\tau\Phi)\le\Psi.
\end{equation}

\smallskip
\textbf{Step 2.} \textit{The Schur estimate.}
For $T=\sqrt{b_\tau}R_{\lambda,D}\sqrt{a_\tau}$, set 
$h=\sqrt{a_\tau}\Psi$ and $j=\sqrt{b_\tau}\Phi$.
These weights are continuous and strictly positive on $\overline D$.
Self-adjointness of the Dirichlet resolvent and
\eqref{ClaimD} give
\[
 Th=\sqrt{b_\tau}R_{\lambda,D}(a_\tau\Psi)\le\frac jp,
 \qquad
 T^*j=\sqrt{a_\tau}R_{\lambda,D}(b_\tau\Phi)\le\frac hq.
\]
Proposition~\ref{lem:schur}, with constants $A=1/p$ and $B=1/q$,
therefore gives \eqref{eq:finite-height-bound}.
\end{proof}

The norm bound \eqref{eq:finite-height-bound} is independent of the height and the ball.
We first let the starting height $\tau\to\infty$ and then the radius
$R\to\infty$ to obtain the global bound.

\begin{proof}[Proof of Proposition~\ref{prop:system-resolvent}]
Fix $\lambda>0$ and recall the kernel in \eqref{eq:weighted-resolvent-kernel}.

\smallskip
\textbf{Step 1.} \textit{The limit in the height.}
 For a ball $D\subset\R^d$, we set 
$T_{\lambda,D,\tau}=\sqrt{b_\tau}R_{\lambda,D}\sqrt{a_\tau}$.
Its kernel is
\[
 K_\tau(x,z)=\sqrt{b_\tau(x)}G_{\lambda,D}(x,z)\sqrt{a_\tau(z)}.
\]
Since $u_t>0$, $v_t>0$ in $\Rp$ and $p,q\ge1$,
$a_\tau$ and $b_\tau$ 
are nondecreasing in $\tau$.
Hence, for a.e.~$(x,z)\in D\times D$, $K_\tau(x,z)$ increases
as $\tau\to\infty$ to
\[
 K_D(x,z)=\sqrt{b(x)}G_{\lambda,D}(x,z)\sqrt{a(z)}.
\]
For $0\le f\in L^2(D)$ and a.e.~$x\in D$, letting $\tau\to\infty$
and applying monotone convergence, we obtain
\[
 T_{\lambda,D,\tau}f(x)\uparrow
 H_D(x):=\int_DK_D(x,z)f(z)\dd z.
\]
Letting $\tau\to\infty$, applying Fatou's lemma and using Lemma~\ref{lem:finite-height}, we obtain
\be{IneqHD}
 \|H_D\|_{L^2(D)}^2
 \le\liminf_{\tau\to\infty}\|T_{\lambda,D,\tau}f\|_{L^2(D)}^2
 \le\frac1{pq}\|f\|_{L^2(D)}^2.
\ee
In particular, the limiting integral is finite a.e.~in $D$.

\smallskip
\textbf{Step 2.} \textit{Exhaustion of the tangential space and conclusion.}
For each $R>0$ we extend the Dirichlet kernel on $B_R^d$ by zero outside the ball.
By Proposition~\ref{rem:dirichlet-resolvent} we have 
$K_{B_R^d}\uparrow K_\lambda$ a.e.~in $\R^d\times\R^d$ as $R\to\infty$.
For a nonnegative compactly supported $f\in L^2(\R^d)$,
letting $R\to\infty$ in \eqref{IneqHD} with $D=B_R^d$ and applying monotone convergence and Fatou's
lemma, we obtain
\[
 \left\|\int_{\R^d}K_\lambda(\cdot,z)f(z)\dd z\right\|_2
 \le(pq)^{-1/2}\|f\|_2.
\]
By standard density arguments, it follows that $K_\lambda$ defines a bounded operator $T_\lambda$ on $L^2(\R^d)$
satisfying \eqref{eq:system-resolvent-bound}. 
\end{proof}

\subsection{Completion of proof of Theorem~\ref{thm:system}}
\label{sec:main-proofs}

Assume for contradiction that $(u,v)\not\equiv(0,0)$ satisfies the hypotheses of Theorem~\ref{thm:system}.
By Remark~\ref{rem:linear} we may assume $pq>1$.
Let $(U,V)$ be the normal limit given by Proposition~\ref{prop:system-limits}.
Set $\bu=(u,v)^{\mathsf T}$, $\bU=(U,V)^{\mathsf T}$ and
\begin{equation}\label{eq:system-ab}
 (a(y),b(y))=(V^{p-1}(y),U^{q-1}(y)),\qquad
  \bC(y)=\begin{pmatrix}0&a(y)\\ b(y)&0\end{pmatrix}.
\end{equation}
The limiting equations give
$\bC\bU=(V^p,U^q)^{\mathsf T}\in L^1_{\loc}(\R^d)^2$.
For every $t>0$,
$u_t>0$ and $v_t>0$ in $\Rp$ give
$0\le u(y,t)\le U(y)$ and $0\le v(y,t)\le V(y)$ for a.e.~$y\in\R^d$. Since
$p,q\ge1$, we have
 $v^p\le V^{p-1}v=av$ and $u^q\le U^{q-1}u=bu$.
Consequently, \eqref{eq:common-hypotheses} holds, and it
follows from Proposition~\ref{lem:common-projection} with $m=2$
that, for every $T>0$, the projections $F_{T,1},F_{T,2}$ defined
in \eqref{eq:common-positive-parts} satisfy
\eqref{eq:common-projection} and \eqref{eq:common-domination}
with $\lambda_T=\frac{\pi^2}{T^2}$. 
In particular 
$$bF_{T,1}^2\le c_TU^{q+1},\quad aF_{T,2}^2\le  c_TV^{p+1}
\quad\hbox{with $c_T=\frac{4T^2}{\pi^2}$.}$$
We will apply Proposition~\ref{lem:system-extinction}
with $(F_1,F_2,\bar a,\bar b)=(F_{T,1},F_{T,2},a,b)$ and $\lambda=\lambda_T$. 
Since $U^{q+1}, V^{p+1}\in\Pclass$ by \eqref{eq:system-limit-energy} in Proposition~\ref{prop:system-limits},
we have $U, V, U^q, V^p \in\Pclass$ by H\"older's inequality
and assumption \eqref{newcondaF} is satisfied. 
On the other hand, \eqref{eq:system-resolvent-bound} in Proposition~\ref{prop:system-resolvent} 
guarantees assumption \eqref{hypTmuL2} with $\theta=(pq)^{-1/2}\in(0,1)$ for every $\mu>0$.
We may thus apply Proposition~\ref{lem:system-extinction} to deduce, for every $T>0$,
that $F_{T,1}=F_{T,2}=0$ a.e.~in $\R^d$, hence 
$(u,v)\equiv(0,0)$ in $\Rp$ by Lemma~\ref{lem:common-finish}.
This contradicts $(u,v)\not\equiv(0,0)$ and concludes the proof of Theorem~\ref{thm:system}.
\qed

\appendix 
\section{Proof of Lemma~\ref{lem:kato} and Propositions~\ref{lem:resolvent-P}--\ref{lem:schur}}

For the convenience of readers, we here provide the proofs of several essentially known results that we used. 

\begin{proof}[Proof of Lemma~\ref{lem:kato}]
Fix $0\le\varphi\in C_c^\infty(\Omega)$ and choose an open set
$\omega\Subset\Omega$ containing its support. Let $\rho_j$ be a
nonnegative standard mollifying sequence. For $j$ large,
$w_j=w*\rho_j$ and $g_j=g*\rho_j$ are defined on $\omega$, and
$\Delta w_j\ge g_j$ there.

Let $\beta\in C^2(\R)$ satisfy $\beta(0)=0$,
$\beta',\beta''\ge0$, and $\beta'\in L^\infty(\R)$. Then
\[
 \Delta\beta(w_j)
 =\beta'(w_j)\Delta w_j+\beta''(w_j)|\nabla w_j|^2
 \ge\beta'(w_j)g_j.
\]
Testing with $\varphi$ and letting $j\to\infty$ along a subsequence
with $w_j\to w$ a.e.~in $\omega$, we obtain
\be{KatoBeta}
 \int_\Omega\beta(w)\Delta\varphi
 \ge\int_\Omega\beta'(w)g\varphi.
\ee
Indeed, the left-hand side converges by the Lipschitz continuity
of $\beta$ and local $L^1$ convergence as $j\to\infty$. On the right, the
difference term arising from replacing $g_j$ by $g$ is bounded by
$\|\beta'\|_\infty\|\varphi\|_\infty\|g_j-g\|_{L^1(\omega)}$,
and the remaining term converges as $j\to\infty$ by dominated convergence.

Choose a nondecreasing $\vartheta\in C^1(\R)$ such that
$\vartheta=0$ on $(-\infty,0]$, $\vartheta=1$ on $[1,\infty)$,
and $0\le\vartheta\le1$, and put
$\beta_k(r)=\int_0^r\vartheta(ks)\dd s$. Then
$0\le\beta_k(r)\le r_+$,
$\beta_k(r)\to r_+$, and
$\beta_k'(r)\to\ind_{\{r>0\}}$ pointwise as $k\to\infty$.
Applying inequality \eqref{KatoBeta}
to $\beta_k$, letting $k\to\infty$
and using dominated convergence, we obtain
$\int_\Omega w_+\Delta\varphi
 \ge\int_\Omega\ind_{\{w>0\}}g\varphi$,
which is \eqref{eq:kato}.
\end{proof}

\begin{proof}[Proof of Proposition~\ref{lem:resolvent-P}]
For $f\in  \Pclass$, we first prove that it defines a tempered distribution.
Let $\int_{B_R^d}|f|\le C(1+R)^N$, take $M>N$, and put
$p_M(\varphi)=\sup_y(1+|y|)^M|\varphi(y)|$. Decomposing into
$B_2^d$ and the dyadic annuli
$A_j=B_{2^{j+1}}^d\setminus B_{2^j}^d$, $j\ge1$, gives
\[
 \int_{\R^d}|f\varphi|\dd y
 \le C p_M(\varphi)
       \Big(1+\sum_{j\ge1}2^{(j+1)N-jM}\Big)
 \le C_{f,M}p_M(\varphi).
\]
The series converges, so $f\in\Sch'(\R^d)$.

We next prove the exponential decay of the kernel away from
its singularity. For $|z|\ge1$, the heat-kernel representation
and $\lambda r/2+|z|^2/(8r)\ge\sqrt\lambda|z|/2$ give
\begin{equation}\label{eq:bessel-exponential-decay}
 G_\lambda(z)
 \le e^{-\sqrt\lambda|z|/2}
       \int_0^\infty(4\pi r)^{-d/2}
                   e^{-\lambda r/2-1/(8r)}\dd r
 \le C_{d,\lambda}e^{-\sqrt\lambda|z|/2}.
\end{equation}
For $R\ge1$, Tonelli's theorem yields
\[
 \int_{B_R^d}(G_\lambda*|f|)(y)\dd y
 =\int_{\R^d}|f(z)|\Big(\int_{B_R^d}G_\lambda(y-z)\dd y\Big)\dd z.
\]
For $z\in\R^d$, we have
$\int_{B_R^d}G_\lambda(y-z)\dd y
 \le\int_{\R^d}G_\lambda(y-z)\dd y=\frac1\lambda$.
Consequently, we obtain
\[
 \int_{|z|\le2R}|f(z)|
       \Big(\int_{B_R^d}G_\lambda(y-z)\dd y\Big)\dd z
       \le\frac1\lambda\int_{B_{2R}^d}|f(z)|\dd z
 \le C(1+R)^N.
\]
For $y\in B_R^d$ and $|z|>2R$,
we have $|y-z|\ge |z|-|y|\ge\frac{|z|}{2}\ge1$.
Using \eqref{eq:bessel-exponential-decay}, we obtain
\[
 \begin{aligned}
 \int_{|z|>2R}|f(z)|
       \Big(\int_{B_R^d}G_\lambda(y-z)\dd y\Big)\dd z
 &\le CR^d\int_{|z|>2R}|f(z)|
                         e^{-\sqrt\lambda|z|/4}\dd z\\
 &\le CR^d\sum_{j\ge1}e^{-\sqrt\lambda j/4}
                    \int_{j\le|z|<j+1}|f(z)|\dd z\\
 &\le CR^d\sum_{j\ge1}(j+2)^N e^{-\sqrt\lambda j/4}
 \le CR^d.
 \end{aligned}
\]
Consequently, $G_\lambda*f\in \Pclass$ is locally integrable. 
The same calculation justifies Fubini's theorem
against compactly supported smooth test functions. 
For such a test $\varphi$, symmetry of $G_\lambda$ gives
\[
 \langle G_\lambda*f,\varphi\rangle
 =\langle f,R_\lambda\varphi\rangle
 =\langle R_\lambda f,\varphi\rangle.
\]
For every $\varphi\in\Sch(\R^d)$, we have
$R_\lambda\varphi\in\Sch(\R^d)$, so $\langle f,R_\lambda\varphi\rangle$
is well defined. For $f\in\Sch'(\R^d)$, the definition of $R_\lambda f$ gives
$\langle R_\lambda f,\varphi\rangle=\langle f,R_\lambda\varphi\rangle$.
Hence $G_\lambda*f=R_\lambda f$ in $\D'(\R^d)$.
The multiplier identity
$(|\xi|^2+\lambda)\widehat{R_\lambda f}=\widehat f$ proves the
distributional equation and proves \eqref{eq:resolvent-P}, completing the proof.
\end{proof}

\begin{proof}[Proof of Proposition~\ref{rem:dirichlet-resolvent}]
Set $G_{\lambda,R}=G_{\lambda,B_R^d}$. 
To prove monotonicity, fix $z\in B_{R_1}^d$ and $0<R_1<R_2$.
The difference
$w=G_{\lambda,R_2}(\cdot,z)-G_{\lambda,R_1}(\cdot,z)$ 
satisfies $(-\Delta_y+\lambda)w=0$ in $B_{R_1}^d$, since the two
Dirac masses cancel. Interior elliptic regularity makes this difference
smooth across~$z$. Its boundary values are
$G_{\lambda,R_2}(\cdot,z)\ge0$ on $\partial B_{R_1}^d$.
For $x\in B_{R_1}^d$ with $x\ne z$, applying the maximum principle,
we obtain $G_{\lambda,R_1}(x,z)\le G_{\lambda,R_2}(x,z)$.

For convergence, fix $x,z\in\R^d$, $x\ne z$, and
$R>\max\{|x|,|z|\}$. Set
$H_R(y)=G_\lambda(y-z)-G_{\lambda,R}(y,z)$.
Again the Dirac masses cancel, so $H_R$ extends smoothly across $z$,
and we have
\[
 \begin{cases}
 (-\Delta_y+\lambda)H_R=0&\text{in }B_R^d,\\
 H_R(\xi)=G_\lambda(\xi-z)\ge0&\text{on }\partial B_R^d.
 \end{cases}
\]
Let $M_R(z)=\sup_{\xi\in\partial B_R^d}G_\lambda(\xi-z)$.
We have $(-\Delta_y+\lambda)(M_R(z)-H_R)=\lambda M_R(z)\ge0$
in $B_R^d$ and $M_R(z)-H_R\ge0$ on $\partial B_R^d$.
Applying the maximum principle to $H_R$ and $M_R(z)-H_R$, we obtain
$
 0\le H_R(x)\le M_R(z).
$
For $R>|z|+1$, we have $|\xi-z|\ge R-|z|>1$ on $\partial B_R^d$.
Using~\eqref{eq:bessel-exponential-decay}, we obtain
$M_R(z)\le C_{d,\lambda}e^{-\sqrt\lambda(R-|z|)/2}\to 0$
 as $R\to\infty$.
This proves \eqref{eq:dirichlet-exhaustion} for every $d\ge1$.
\end{proof}

\begin{proof}[Proof of Proposition~\ref{lem:resolvent-upper-bound}]
The proof is done by a cutoff iteration, using polynomial integral growth. 
Applying Proposition~\ref{lem:resolvent-P}, we have
$H:=F-R_\lambda f\in\Pclass$. Subtracting the equation for
$R_\lambda f$ from the inequality for $F$, we obtain
$(-\Delta_y+\lambda)H\le0$ in $\D'(\R^d)$. Kato's inequality (Lemma~\ref{lem:kato}) gives
$\Delta_yH_+\ge\lambda H_+$.
Choose $\eta_R\in C_c^\infty(B_{2R}^d)$ with
$0\le\eta_R\le1$, $\eta_R=1$ on $B_R^d$, and
$|\Delta\eta_R|\le c/R^2$. Testing this distributional
inequality with $\eta_R$ yields
\[
 \lambda\int_{B_R^d}H_+\dd y
 \le\lambda\int_{\R^d}H_+\eta_R\dd y
 \le\int_{\R^d}H_+\Delta\eta_R\dd y
 \le\frac c{R^2}\int_{B_{2R}^d}H_+\dd y.
\]
Writing $M(R)=\int_{B_R^d}H_+\dd y$ and iterating, we obtain,
for every integer $k\ge1$,
\[
 M(R)\le\left(\frac c{\lambda R^2}\right)^k
             2^{-k(k-1)}M(2^kR).
\]
For each fixed $R\ge1$, polynomial integral growth gives
$M(2^kR)\le C_R2^{kN}$. The right-hand side therefore tends
to zero as $k\to\infty$, so $H_+=0$ a.e.~in $\R^d$.
This proves \eqref{eq:resolvent-upper-bound}. 
\end{proof}

\begin{proof}[Proof of Proposition~\ref{lem:schur}]
Fix $f\in L^2(D)$ and define the nonnegative, possibly infinite integral
\[
 H_f(x)=\int_D K(x,z)|f(z)|\dd z.
\]
For integers $m\ge1$, set
$H_{f,m}(x)=\int_D\min\{K(x,z),m\}\min\{|f(z)|,m\}\dd z$.
As $m\to\infty$, we have $H_{f,m}\uparrow H_f$.
Applying Cauchy--Schwarz  inequality to $H_{f,m}$, letting $m\to\infty$
and using monotone convergence, we obtain
\[
 H_f(x)^2
 \le A j(x)\int_D K(x,z)\frac{|f(z)|^2}{h(z)}\dd z.
\]
Using Tonelli's theorem and the second inequality in
\eqref{eq:schur-hypotheses}, we obtain
\[
 \int_D H_f^2\dd x
\le A\int_D\frac{|f(z)|^2}{h(z)}
          \left(\int_D K(x,z)j(x)\dd x\right)\dd z
          \le AB\int_D|f(z)|^2\dd z.
\]
Thus $H_f\in L^2(D)$ and $H_f<\infty$ a.e.~in $D$. The integral
$Tf(x)=\int_D K(x,z)f(z)\dd z$ is therefore absolutely convergent
for a.e.~$x\in D$. It defines a linear operator, and
$|Tf|\le H_f$ gives \eqref{eq:schur-conclusion}.
\end{proof}

{\bf Acknowledgements.} 
Yimei Li is supported by the Fundamental Research Funds for the Central Universities (No.~2026JBMC040) and the National Natural Science Foundation of China (No.~12371200).

\smallskip

{\bf Declarations.} The authors report that there are no competing interests to declare. This manuscript has no associated data.

\smallskip

{\bf Declaration on the use of generative AI.}
The authors used ChatGPT to assist them in literature review, language editing and manuscript presentation.
They take full responsibility for the manuscript.

\end{document}